\documentclass[english]{amsart}
\usepackage[T1]{fontenc}
\usepackage[latin9]{inputenc}
\usepackage{float}
\usepackage{units}
\usepackage{amsthm}
\usepackage{amsbsy}
\usepackage{graphicx}
\usepackage{amssymb}
\usepackage[all]{xy}

\numberwithin{equation}{section} 
\numberwithin{figure}{section} 
\theoremstyle{plain}
\theoremstyle{plain}
\newtheorem{thm}{Theorem}
  \theoremstyle{definition}
  \newtheorem{defn}[thm]{Definition}
  \theoremstyle{plain}
  \newtheorem{prop}[thm]{Proposition}
  \theoremstyle{remark}
  \newtheorem*{rem*}{Remark}
  \theoremstyle{plain}
  \newtheorem{cor}[thm]{Corollary}

\usepackage{nicefrac}

\date{\today}

\def\at{\Big|}
\def\rv{\textrm{rv}}
\def\ind{{\rm Ind}}

\usepackage[all]{xy}

\usepackage{babel}

\begin{document}

\title{Linear Representations and Isospectrality with Boundary Conditions}

\author{Ori Parzanchevski$^{1}$ and Ram Band$^{2}$}

\address{$^{\text{1}}$ Institute of Mathematics\\
Hebrew University, Jerusalem 91904, Israel}

\email{parzan@math.huji.ac.il}

\address{$^{\text{2}}$ Department of Physics of Complex Systems\\
The Weizmann Institute of Science, Rehovot 76100, Israel}

\email{rami.band@weizmann.ac.il}
\begin{abstract}
We present a method for constructing families of isospectral systems,
using linear representations of finite groups. We focus on quantum
graphs, for which we give a complete treatment. However, the method
presented can be applied to other systems such as manifolds and two-dimensional
drums. This is demonstrated by reproducing some known isospectral
drums, and new examples are obtained as well. In particular, Sunada's
method \cite{Sunada} is a special case of the one presented.
\end{abstract}
\maketitle

\section{Introduction}

\label{sec:introduction}{}``Can one hear the shape of a drum?''
- This question was posed by Marc Kac in 1966 \cite{Kac}. In other
words, is it possible to determine the shape of a planar Euclidean
domain from the spectrum of the Laplace operator on it? This question
gave rise to fertile research, investigating it from various aspects.
Two main approaches were, on the one hand, attempts to deal with the
inverse question of reconstructing the shape from the spectrum, and
on the other hand, trying to find systems whose shapes are different,
yet have the same spectrum. Such examples are called isospectral.
Although Kac's original question regarded two dimensional planar drums,
the research on isospectrality expanded quickly to other types of
systems. We will not go into detail, but refer the interested reader
to \cite{Sunada}-\cite{Brooks-Sun} for a broader view of the field.
However, we will mention here two milestones in the field of isospectrality.
A theorem by Sunada gave an important machinery for the construction
of isospectral Riemannian manifolds \cite{Sunada}. Later, this method
was used by Gordon, Webb and Wolpert to construct the first pair of
isospectral planar Euclidean domains \cite{Gordon1,Gordon2} thus
negatively answering Kac's original question.

This paper starts with a presentation of the basic theory of quantum
graphs and existing results on quantum graph isospectrality. We then
present the algebraic part of our theory and its main theorem. This
is followed by a section which explains the construction of the so
called quotient graphs that lie in the heart of the theory. After
the theory is fully presented, we apply it to obtain various examples
of isospectral quantum graphs. We then demonstrate how to apply the
method to other systems, explaining some known results, as well as
obtaining new ones. In particular we discuss the relation to Sunada's
method. We conclude by pointing out key elements of the theory that
are to be investigated further and by presenting open questions.

\section{Quantum graphs}

\label{sec:graphs}A \emph{graph} $\Gamma$ consists of a finite set
of vertices $V=\left\{ v_{i}\right\} $ and a finite set $E=\left\{ e_{j}\right\} $
of edges connecting the vertices. We assume that there are no parallel
edges (different edges with the same endpoints) or loops (edges connecting
a vertex to itself), but we shall see that this inflicts only a small
loss of generality. We denote by $E_{v}$ the set of all edges incident
to the vertex $v$. The degree (valency) of the vertex $v$ is $d_{v}=\left|E_{v}\right|$.
$\Gamma$ becomes a \emph{metric graph} if each edge $e\in E$ is
assigned a finite length $l_{e}>0$. It is then possible to identify
the edge $e$ with a finite segment $\left[0,l_{e}\right]$ of the
real line, having the natural coordinate $x_{e}$ along it. A function
on the graph is a vector $f=\left(f\big|_{e_{1}},\ldots,f\big|_{e_{\left|E\right|}}\right)$
of functions $f\big|_{e_{j}}:\left[0,l_{e_{j}}\right]\rightarrow\mathbb{C}$
on the edges. We shall usually consider smooth functions on the graph,
meaning that $f\big|_{e}\in C^{^{\infty}}\negthickspace\left(\left[0,l_{e}\right]\right)$
for all $e\in E$. Notice that in general it is not required that
for $v\in V$ and $e,e'\in E_{v}$ the functions $f\big|_{e}$ and
$f\big|_{e'}$ agree on $v$.

To obtain a \emph{quantum graph}, we consider a differential operator
on the graph, by default the negative Laplacian: $-\Delta f=\left(-f''\big|_{e_{1}},\ldots,-f''\big|_{e_{\left|E\right|}}\right)$.
In addition, we require the functions on the graph to obey certain
boundary conditions stated a priori; for each vertex $v\in V$, we
consider homogeneous boundary conditions which involve the values
and derivatives of the function at the vertex, of the form $A_{v}\cdot f\big|_{v}+B_{v}\cdot f'\big|_{v}=0$.
Here $A_{v}$ and $B_{v}$ are $d_{v}\times d_{v}$ complex matrices,
$f\big|_{v}$ is the vector $\left(\begin{matrix}f\big|_{e_{n_{1}}}(v) & \ldots & f\big|_{e_{n_{d_{v}}}}(v)\end{matrix}\right)^{T}$
of the values of $f$ on the edges in $E_{v}$ at $v$, and $f'\big|_{v}=\left(\begin{matrix}f'\big|_{e_{n_{1}}}(v) & \ldots & f'\big|_{e_{n_{d_{v}}}}(v)\end{matrix}\right)^{T}$
is the vector of outgoing derivatives of $f$ taken at the vertex.
To sum up, a quantum graph is a metric graph equipped with a differential
operator and homogeneous differential boundary conditions at the vertices.
Notice that before stating the boundary conditions, the graph is merely
a collection of independent edges with functions defined separately
on each edge. The connectivity of the graph is manifested through
the boundary conditions. We denote by $C^{^{\infty}}\negthickspace\left(\Gamma\right)$
the space of (smooth complex) functions on the graph which satisfy
the boundary conditions at the vertices, and by $\mathcal{H}\left(\Gamma\right)$
the subspace of $C^{^{\infty}}\negthickspace\left(\Gamma\right)$
spanned by eigenfunctions of the Laplacian. The reader interested
in more information about quantum graphs is referred to the reviews
\cite{GS06,Kuch04,Kuch08}.

A standard choice of boundary conditions which we adopt is the so
called \textit{Neumann} boundary condition%
\footnote{This condition is also widely encountered under the name of \emph{Kirchhoff}
condition.%
}:\\
 $\bullet\quad$ $f$ agrees on the vertices: $\forall v\in V\quad\forall e,e'\in E_{v}\::\quad f\at_{e}\left(v\right)=f\at_{e'}\left(v\right)$.
\\
 $\bullet\quad$ The sum of outgoing derivatives at each vertex
is zero: $\forall v\in V$ :$\sum\limits _{e\in E_{v}}f'\at_{e}\left(v\right)=0$.\\
The Neumann boundary condition can thus be represented by the matrices\[
A_{v}=\left(\begin{matrix}1 & -1\\
 & \ddots & \ddots\\
 &  & 1 & -1\\
0 & \cdots & 0 & 0\end{matrix}\right)\,,\qquad B_{v}=\left(\begin{matrix}0 & 0 & \cdots & 0\\
\vdots & \vdots & \ddots & \vdots\\
0 & 0 & \cdots & 0\\
1 & 1 & \cdots & 1\end{matrix}\right)\:.\]
For a vertex of degree one the Neumann condition is expressed by the
matrices $A_{v}=\left(\begin{matrix}0\end{matrix}\right)$, $B_{v}=\left(\begin{matrix}1\end{matrix}\right)$,
and means that the derivative of the function is zero at the leaf
$v$. Another natural boundary condition for leaves is the \emph{Dirichlet
}boundary condition: $A_{v}=\left(\begin{matrix}1\end{matrix}\right)$,
$B_{v}=\left(\begin{matrix}0\end{matrix}\right)$, which means that
the function vanishes at the vertex.

Neumann vertices of degree two deserve a special attention. They can
be thought of as inner points along a single edge - the concatenation
of the two edges incident to the vertex - and we would like to be
able to add or remove such inner points, for reasons which will become
clear later on. At such points, however, a function on the edge is
only required (by the Neumann condition) to be continuously differentiable
($C^{1}$), rather then smooth ($C^{^{\infty}}$); therefore, adding
a Neumann vertex of degree two at an inner point of an edge augments
the space of allowed functions (by ones such as $\left|x\right|\cdot x$).
The question of $C^{1}$ versus $C^{^{\infty}}$ is inherent to the
modeling of one dimensional manifolds as quantum graphs. For example,
in order to regard the circle $S^{1}$ as a quantum graph, we must
place at least one vertex along it, and at this vertex functions on
the resulting graph may have a non-differentiable derivative.

The good news is that adding or removing Neumann vertices of degree
two does not change the \textbf{spectral }properties of the graph
in question. For sums of eigenfunctions of the Laplacian, being $C^{1}$
and piecewise $C^{^{\infty}}$ is equivalent to being $C^{^{\infty}}$
altogether, so that if the graph $\Gamma'$ is obtained from $\Gamma$
by adding or removing such points, we have $\mathcal{H}\left(\Gamma'\right)=\mathcal{H}\left(\Gamma\right)$.
With this observation in mind, we will allow ourselves to make manipulations
of degree two Neumann vertices, with no essential loss of generality
(at least from the spectral viewpoint). For example, loops and parallel
edges can be eliminated by the introduction of such {}``dummy''
vertices, so that as mentioned, we shall assume that we are dealing
with graphs with no such nuisances.

If for every $v\in V$ the $d_{v}\times2d_{v}$ matrix $\left(A_{v}\,\middle|\, B_{v}\right)$
is of full rank, we shall say that the quantum graph is \emph{exact}.
Non-exact quantum graphs are not very interesting from the spectral
point of view, as their spectrum is all of $\mathbb{C}$. On the other
hand, we shall later be led to consider the opposite phenomena, i.e.,
vertices at which there are {}``too many'' boundary conditions.
In this case we shall admit $A_{v}$ and $B_{v}$ to be of size $m\times d_{v}$,
possibly with $m>d_{v}$, and we shall call the corresponding graphs
\emph{generalized quantum graphs}. From the spectral perspective these
are much more interesting than non-exact quantum graphs. Consider
for example a Y-shaped graph, with a Neumann condition at the center,
Dirichlet conditions at two of the leaves, and the condition $A_{v}=\left(\begin{smallmatrix}1\\
0\end{smallmatrix}\right)$, $B_{v}=\left(\begin{smallmatrix}0\\
1\end{smallmatrix}\right)$ at the third; its spectrum is nonempty if and only if the lengths
of the two edges with Dirichlet leaves are commensurable.

There is a natural inner product on $C^{^{\infty}}\negthickspace\left(\Gamma\right)$,
given by $\left\langle f,g\right\rangle =\sum\limits _{e\in E}\int_{0}^{l_{e}}\overline{f\big|_{e}}\cdot g\big|_{e}dx_{e}$.
Kostrykin and Schrader \cite{KosSch99} provide necessary and sufficient
conditions for the Laplacian to be self-adjoint with respect to this
product. These conditions can be stated in a number of equivalent
forms (see \cite{Kuch08}). We give two of them:
\begin{enumerate}
\item $\Gamma$ is exact, and $A_{v}\cdot B_{v}^{\dagger}$ is self-adjoint
for every $v\in V$.
\item For every $v\in V$ there exist a unitary matrix $U$ such that $\left(A_{v}\,\middle|\, B_{v}\right)$
is row-equivalent to $\left(i\left(U-I\right)\,\middle|\, U+I\right)$
\footnote{(2) follows from (1) by taking $U=2\left(-A+iB\right)^{-1}A+I$. The
other direction is trivial.%
}.
\end{enumerate}
In particular, Neumann and Dirichlet boundary conditions satisfy these
requirements.

There are several known results concerning isospectrality of quantum
graphs. Gutkin and Smilansky \cite{gutkinus} show that under certain
conditions a quantum graph can be heard, meaning that it can be recovered
from the spectrum of its Laplacian. On the other hand, constructions
of isospectral graphs were also established, by various means: by
a trace formula for the heat kernel \cite{Roth}, by turning isospectral
discrete graphs into equilateral quantum graphs \cite{vonBelow},
and weighted discrete graphs into non-equilateral ones \cite{Oren};
in \cite{gutkinus,Tashkent} a wealth of examples is given by an analogy
to the isospectral drums obtained by Buser et al.\ \cite{BuserConway},
and in \cite{BSS06} is presented an example, whose generalization
has led to the theory presented in this paper.

\section{Algebra}

\label{sec:algebra}For a quantum graph $\Gamma$ we can regard $\mathcal{H}\left(\Gamma\right)$
as a $\mathbb{C}\left[x\right]$-module, where $x$ acts as the (negative)
Laplacian. For every $\lambda\in\mathbb{C}$ we denote by $\Phi_{\Gamma}\left(\lambda\right)$
the submodule \[
\Phi_{\Gamma}\left(\lambda\right)=\mathrm{Ann}_{\mathcal{H}\left(\Gamma\right)}\left(x-\lambda\right)=\left\{ f\in\mathcal{H}\left(\Gamma\right)|-\Delta f=\lambda f\right\} \,,\]
which as a vector space is merely the $\lambda$-eigenspace of the
Laplacian. \emph{The spectrum of $\Gamma$} is the function \[
\sigma_{\Gamma}:\lambda\mapsto\dim_{\mathbb{C}}\Phi_{\Gamma}\left(\lambda\right)\,,\]
which assigns to each eigenvalue its multiplicity%
\footnote{In effect we have $\sigma_{\Gamma}:\mathbb{C}\rightarrow\left\{ 0..2\left|E\right|\right\} $,
as the eigenvalue of a Laplacian eigenfunction, together with the
values $\left\{ f\big|_{e}\left(0\right),\, f'\big|_{e}\left(0\right)\right\} _{e\in E}$,
determine the function.%
}. Two quantum graphs $\Gamma$ and $\Gamma'$ are said to be \emph{isospectral}
if their spectra coincide, that is $\sigma_{\Gamma}\equiv\sigma_{\Gamma'}$,
and as noted in \cite{Webster}, this can follow from the stronger
assumption that $\mathcal{H}\left(\Gamma\right)$ and $\mathcal{H}\left(\Gamma'\right)$
are isomorphic as $\mathbb{C}\left[x\right]$-modules, which means
that their Laplacians are conjugate.

A symmetry of a quantum graph is an invertible graph map that preserves
both the lengths of edges and the boundary conditions at the vertices.
The group of all such symmetries is denoted $\mathrm{Aut}\,\Gamma$.
A left action of a group $G$ on a quantum graph $\Gamma$ is equivalent
to a group homomorphism $G\rightarrow\mathrm{Aut}\,\Gamma$. Such
action induces a left action of $G$ on $\mathcal{H}\left(\Gamma\right)$
(by $\left(gf\right)\left(x\right)=f\left(g^{-1}x\right)$ - the inversion
accounts for the contravariantness of $\mathcal{H}$). This gives
$\mathcal{H}\left(\Gamma\right)$ a $\mathbb{C}G\left[x\right]$-module
structure, since the Laplacian commutes with all symmetries. The eigenspaces
$\Phi_{\Gamma}\left(\lambda\right)=\mathrm{Ann}\left(x-\lambda\right)$
are again submodules, and in particular they are $\mathbb{C}G$-modules,
that is, complex representations of $G$. Assuming that $G$ is finite,
with irreducible complex representations $S_{1},\ldots,S_{r}$, we
can decompose each eigenspace to its isotypic components: \begin{equation}
\Phi_{\Gamma}\left(\lambda\right)=\bigoplus\limits _{i=1}^{r}\Phi_{\Gamma}^{S_{i}}\left(\lambda\right)\,,\label{eq:isotypic-decomp}\end{equation}
 where $\Phi_{\Gamma}^{S_{i}}\left(\lambda\right)\cong S_{i}\oplus\ldots\oplus S_{i}$
as $\mathbb{C}G$-modules.

We start by counting separately, for each irreducible representation
$S$ of $G$, only the $\lambda$-eigenfunctions which reside in $\Phi_{\Gamma}^{S}\left(\lambda\right)$.
This means that we are restricting our attention to functions which
under the action of $\mathbb{C}G$ span a space that is isomorphic,
as a representation of $G$, to $S$. However, since $\dim S$ always
divides $\dim\Phi_{\Gamma}^{S}\left(\lambda\right)$, we can already
normalize by it. We thus define \emph{the spectrum of $S$} as \begin{equation}
\sigma_{\Gamma}^{S}:\lambda\mapsto\nicefrac{\dim_{\mathbb{C}}\Phi_{\Gamma}^{S}\left(\lambda\right)}{\dim_{\mathbb{C}}S}\:.\label{eq:spectrum-of-S}\end{equation}
 By the orthogonality relations of irreducible characters, we can
rewrite this as $\sigma_{\Gamma}^{S}\left(\lambda\right)=\left\langle \chi_{S},\chi_{\Phi_{\Gamma}\left(\lambda\right)}\right\rangle _{G}$,
and expanding this linearly, we define \emph{the spectrum of $R$},
for every representation $R$ of $G$, to be \begin{equation}
\sigma_{\Gamma}^{R}:\lambda\mapsto\left\langle \chi_{R},\chi_{\Phi_{\Gamma}\left(\lambda\right)}\right\rangle _{G}\:.\label{eq:spectrum-of-R}\end{equation}

$\sigma_{\Gamma}^{S}\left(\lambda\right)$ has an algebraic significance:
it reflects the size of the $S$-isotypic part of $\Phi_{\Gamma}\left(\lambda\right)$.
Looking for a parallel algebraic interpretation of $\sigma_{\Gamma}^{R}\left(\lambda\right)$,
we find that \[
\sigma_{\Gamma}^{R}\left(\lambda\right)=\dim_{\mathbb{C}}\mathrm{Hom}_{\mathbb{C}G}\left(R,\Phi_{\Gamma}\left(\lambda\right)\right)\:.\]

Quite generally, if $A\hookrightarrow B$ is a ring extension and
$M$ and $N$ are modules over $A$ and $B$ respectively, then for
$C_{B}\left(A\right)$, the centralizer of $A$ in $B$~%
\footnote{That is, $C_{B}\left(A\right)=\left\{ b\in B\,\middle|\,\forall a\in A:ab=ba\right\} $.%
}, $\mathrm{Hom}_{A}\left(M,N\right)$ has a natural $C_{B}\left(A\right)$-module
structure (by $\left(bf\right)\left(m\right)=b\cdot f\left(m\right)$
for every $b\in C_{B}\left(A\right)$). For our purposes, since $\mathbb{C}\left[x\right]$
lies in the centralizer of $\mathbb{C}G$ in $\mathbb{C}G\left[x\right]$,
we obtain that $\mathrm{Hom}_{\mathbb{C}G}\left(R,\mathcal{H}\left(\Gamma\right)\right)$
has a $\mathbb{C}\left[x\right]$-module structure:\[
\hphantom{\left(\tilde{f}\in\mathrm{Hom}_{\mathbb{C}G}\left(R,\mathcal{H}\left(\Gamma\right)\right)\right)}x\cdot\tilde{f}:r\mapsto-\Delta\left(\tilde{f}\left(r\right)\right)\qquad\left(\tilde{f}\in\mathrm{Hom}_{\mathbb{C}G}\left(R,\mathcal{H}\left(\Gamma\right)\right)\right)\,.\]
This allows us to make the following definition.
\begin{defn}
\label{def:quotient-definition}A \emph{$\nicefrac{\Gamma}{R}$-graph}
is any quantum graph $\Gamma'$ such that there is a $\mathbb{C}\left[x\right]$-module
isomorphism\begin{equation}
\mathcal{H}\left(\Gamma'\right)\cong\mathrm{Hom}_{\mathbb{C}G}\left(R,\mathcal{H}\left(\Gamma\right)\right)\:.\label{eq:quotient-definition}\end{equation}

\end{defn}
We note, in particular, that for such $\Gamma'$, there is an isomorphism\begin{eqnarray*}
\Phi_{\Gamma'}\left(\lambda\right) & =\\
\mathrm{Ann}_{\mathcal{H}\!\left(\Gamma'\right)}\left(x-\lambda\right) & \cong & \mathrm{Ann}_{\mathrm{Hom}_{\mathbb{C}G}\left(R,\mathcal{H}\!\left(\Gamma\right)\right)}\left(x-\lambda\right)\\
 & = & \mathrm{Hom}_{\mathbb{C}G}\left(R,\Phi_{\Gamma}\left(\lambda\right)\right)\end{eqnarray*}
 which by taking dimensions translates to equality of spectra: \begin{equation}
\sigma_{\Gamma'}\equiv\sigma_{\Gamma}^{R}\:.\label{eq:spectrum of quotient}\end{equation}

Since $\sigma_{\Gamma}^{R}$ is not a spectrum in the classical sense,
we cannot really call this isospectrality. However we do have from
this that all $\nicefrac{\Gamma}{R}$-graphs are isospectral to one
another, and we will use this to speak non-rigorously about {}``the
spectrum of $\nicefrac{\Gamma}{R}$'', $\sigma_{\nicefrac{\Gamma}{R}}\equiv\sigma_{\Gamma}^{R}$.
The following proposition exhibits another manifestation of isospectrality.
\begin{prop}
\label{pro:regular-spec}All $\nicefrac{\Gamma}{\mathbb{C}G}$-graphs
are isospectral to $\Gamma$. \end{prop}
\begin{proof}
By \eqref{eq:isotypic-decomp}, \eqref{eq:spectrum-of-S}, and linearity,
the classical spectrum $\sigma_{\Gamma}$ coincides with the spectrum
of the regular representation of $G$: \begin{equation}
\sigma_{\Gamma}^{\mathbb{C}G}\equiv\sum_{i=1}^{r}\dim S_{i}\cdot\sigma_{\Gamma}^{S_{i}}\equiv\sigma_{\Gamma}\:.\label{eq:regular-spectrum}\end{equation}
This can also be deduced from the fact that for every $R$-module
$M$ there is an isomorphism $\mathrm{Hom}_{R}\left(R,M\right)\cong M$,
so that we have $\mathrm{Hom}_{\mathbb{C}G}\left(\mathbb{C}G,\Phi_{\Gamma}\left(\lambda\right)\right)\cong\Phi_{\Gamma}\left(\lambda\right)$
for every eigenvalue $\lambda$.
\end{proof}
We can say even more:
\begin{thm}
\label{thm:mainthm}Let $\Gamma$ be a quantum graph equipped with
an action of $G$, $H$ a subgroup of $G$, and $R$ a representation
of $H$. Then $\nicefrac{\Gamma}{R}$ is isospectral to $\nicefrac{\Gamma}{\ind_{H}^{G}R}$.\end{thm}
\begin{proof}
This follows at once from the Frobenius Reciprocity Theorem, which
states that there is an isomorphism $\mathrm{Hom}_{\mathbb{C}H}\left(R,\mathcal{H}\left(\Gamma\right)\right)\cong\mathrm{Hom}_{\mathbb{C}G}\left(\ind_{H}^{G}R,\mathcal{H}\left(\Gamma\right)\right)$.
It is straightforward to verify that this is an isomorphism of $\mathbb{C}\left[x\right]$-modules.
Note that from the formal point of view, we have actually shown that
$\nicefrac{\Gamma}{R}$ and $\nicefrac{\Gamma}{\ind_{H}^{G}R}$ are
identical\textbf{ }(as classes\textbf{ }of quantum graphs).\end{proof}
\begin{rem*}
This gives yet another explanation for the equality of the classical
spectrum with that of the regular representation (proposition \ref{pro:regular-spec}):
for $H=\left\{ id\right\} $ and $\mathbf{1}_{H}$ its trivial representation,
it is clear by the isotypic component perspective that $\left(\mathcal{H}\left(\Gamma\right)\right)^{\mathbf{1}_{H}}=\mathcal{H}\left(\Gamma\right)$,
so that \eqref{eq:regular-spectrum} follows from $\ind_{H}^{G}\mathbf{1}_{H}\cong\mathbb{C}G$.\end{rem*}
\begin{cor}
\label{cor:sunadapair}If $G$ acts on $\Gamma$ and $H_{1},\, H_{2}$
are subgroups of $G$ with corresponding representations $R_{1},\, R_{2}$,
such that $\ind_{H_{1}}^{G}R_{1}\cong\ind_{H_{2}}^{G}R_{2}$, then
$\nicefrac{\Gamma}{R_{1}}$ and $\nicefrac{\Gamma}{R_{2}}$ are isospectral. \end{cor}
\begin{rem*}
This corollary is in fact equivalent to the theorem, which follows
by taking $H_{2}=G$, $R_{2}=\ind_{H_{1}}^{G}R_{1}$. It is presented
for being of practical usefulness (it allows one to work with representations
of lower dimension, as can be seen in section \ref{sec:isospectral_qg}),
but also since it indicates the bridge connecting our method with
the classical one of Sunada. In section \ref{sub:Sunada}, we shall
cross it.
\end{rem*}
The sharpest observations in this section would be mere algebraic
tautologies, unless we can show that $\nicefrac{\Gamma}{R}$-graphs
do exist. The next section is devoted to this purpose.

\section{Building $\nicefrac{\Gamma}{R}$-graphs}

\label{sec:quotient graph}In this section we prove the existence
of the quotient graphs $\nicefrac{\Gamma}{R}$. This is done by describing
an explicit construction of $\nicefrac{\Gamma}{R}$, given a graph
$\Gamma$, a representation $R$ of some group $G$ acting on the
graph, and various choices of bases for this representation. As the
lengthy technical details of the construction might encloud the essence
of the method, the reader may prefer to go over section \ref{sec:isospectral_qg}
first, and obtain an intuition for the construction of the quotient
graph from the examples presented there. More intuition for the construction
can be gained from the examples in \cite{BPbS}.

We summarize the main conclusions of this section in the following
theorem:
\begin{thm}
\label{thm:construction-thm}For any representation $R$ of a finite
group $G$, which acts upon a quantum graph $\Gamma$, there exists
a generalized $\nicefrac{\Gamma}{R}$ quantum graph. Furthermore,
if $\Gamma$'s Laplacian is self-adjoint, then there exists a proper
$\nicefrac{\Gamma}{R}$ quantum graph, and it is exact.
\end{thm}

\subsection{\label{sub:construction-intuition}Intuition}

A motivation for the construction of our quotient graphs is given
by thinking about it as an {}``encoding scheme''~%
\footnote{In section \ref{sec:manifolds} we show that the same construction
and motivation can be applied analogously to other geometric systems.%
}. An element $\tilde{f}$ in $\mathrm{Hom}_{\mathbb{C}G}\left(R,\mathcal{H}\left(\Gamma\right)\right)$
can be thought of as a family of functions on $\Gamma$, parametrized
by $R$. To emphasize this view, we shall write $\tilde{f}_{r}$ instead
of $\tilde{f}\left(r\right)$ (where $r\in R$). Our goal is to build
a new graph, each of whose complex functions encodes exactly one such
family. The desired map $\Psi:\mathrm{Hom}_{\mathbb{C}G}\left(R,\mathcal{H}\left(\Gamma\right)\right)\rightarrow\mathcal{H}\left(\nicefrac{\Gamma}{R}\right)$
(see definition \ref{def:quotient-definition}) is in fact this encoding.
An encoding scheme should always be injective (in order to allow decoding),
but we have also required $\Psi$ to be surjective: this can be translated
to the idea that the encoding must be {}``as efficient as possible''~%
\footnote{In a suitable sense, since better encoding may exist, but we want
the encoding to be by another quantum graph, in a manner which intertwines
the corresponding Laplacians.%
} - that $\nicefrac{\Gamma}{R}$ is to be a {}``minimal'' graph allowing
such an encoding, since it admits no Laplacian eigenfunctions apart
from the ones used by the scheme.

First, we reduce the infinite family $\tilde{f}$ to a finite one
by choosing a basis $B=\left\{ b_{j}\right\} _{j=1}^{d}$ for $R$,
and restricting our attention to $\left\{ \tilde{f}_{b_{j}}\right\} _{j=1}^{d}$.
From these {}``basis functions'' we can reconstruct $\tilde{f}$,
since the $\mathbb{C}G$-linearity of $\tilde{f}$ implies in particular
$\mathbb{C}$-linearity (i.e., $\tilde{f}_{\Sigma\alpha_{j}b_{j}}=\sum\alpha_{j}\tilde{f}_{b_{j}}$).
As a first encoding attempt we could take a graph with $d$ times
each edge in $\Gamma$, and let the $j^{{\rm th}}$ copy of the edge
$e$ carry the $j^{{\rm th}}$ basis function restricted to $e$.
That is, define $\left(\Psi\tilde{f}\right)\at_{e_{j}}\equiv\tilde{f}_{b_{j}}\at_{e}$.
However, this encoding is not efficient enough, since we have used
only $\mathbb{C}$-linearity. For each $g\in G$, $\mathbb{C}G$-linearity
implies that $\left\{ \tilde{f}_{r}\at_{e}\right\} _{r\in R}$ determines
$\left\{ \tilde{f}_{r}\at_{ge}\right\} _{r\in R}$, specifically by
\begin{equation}
\tilde{f}_{r}\at_{ge}\equiv\left(g^{-1}\tilde{f}_{r}\right)\at_{e}\equiv\tilde{f}_{g^{-1}r}\at_{e}\label{eq:edge-to-vector}\end{equation}
(the inversion occurs since $G$ acts on $\mathcal{H}\left(\Gamma\right)$
by $g\cdot f=f\circ g^{-1}$). Thus, it suffices to encode the basis
functions on only \textbf{one} edge from each $G$-orbit of $\Gamma$'s
edges.

It turns out that if the action of $G$ on $E$ is free, then apart
from determining the appropriate boundary conditions at the vertices
we are done: for $\left\{ e^{i}\right\} $, a choice of representatives
for $\nicefrac{E}{G}$, setting $\left(\Psi\tilde{f}\right)\at_{e_{j}^{i}}\equiv\tilde{f}_{b_{j}}\at_{e^{i}}$
(where $1\leq j\leq d$) is indeed a {}``good'' encoding (i.e.,
once the boundary conditions are correctly stated, $\Psi$ is an isomorphism.)

If, however, some edge $e=\left\{ v,v'\right\} $ has a non-trivial
stabilizer $G_{e}=G_{v}\cap G_{v'}$, then greater efficiency can
(and therefore must) be achieved. For example, assume that $\dim R=1$
and that for some $g\in G_{e}$ we have $g\notin\ker\rho_{R}$, where
$\rho_{R}$ is the structure homomorphism $G\rightarrow{\rm GL}_{1}\left(\mathbb{C}\right)$.
We then have\[
\tilde{f}_{r}\at_{e}\equiv\tilde{f}_{r}\at_{g^{-1}e}\equiv\tilde{f}_{gr}\at_{e}\equiv\tilde{f}_{\rho_{R}\left(g\right)\cdot r}\at_{e}\equiv\rho_{R}\left(g\right)\cdot\tilde{f}_{r}\at_{e}\]
which implies that $\tilde{f}_{r}\at_{e}\equiv0$ for all $r$, and
as a result, the edge $e$ need not have \textbf{any} representative
in the quotient. We can {}``decode $\tilde{f}_{r}\at_{e}$ from thin
air'', since we know in advance that it can only be the zero function.
The generalization of this observation is that for each edge $e$,
the information in $\left\{ \tilde{f}_{r}\at_{e}\right\} _{r\in R}$
is encapsulated in $R^{G_{e}}$~%
\footnote{$R^{H}$ is the trivial component of ${\rm Res}_{H}^{G}R$, i.e. $R^{H}=\left\{ r\in R\,\middle|\,\forall h\in H\,:\, hr=r\right\} $.%
} : if $r$ belongs to a nontrivial component of ${\rm Res}_{G_{e}}^{G}R$,
then $\tilde{f}_{r}\at_{e}\equiv0$. Therefore, we need only $d_{i}=\dim R^{G_{e^{i}}}$
copies of each representative $e^{i}$ in the quotient%
\footnote{\label{fn:dead-edges}However, we shall later find it convenient to
think about $d=\dim R$ copies, where the $d_{i}+1\ldots d$ copies
are {}``dead'', meaning that whenever a function on them appears
in a formula it is to be understood as zero.%
}. This further {}``compression'' slightly complicates the determination
of the boundary conditions. When $G$ acted freely on the edges, we
had $d$ functions, $\tilde{f}_{b_{j}}$, each satisfying the boundary
conditions at the vertices of $\Gamma$, and we could have translated
this quite easily to boundary conditions on the quotient. Now, however,
for each edge $e^{i}$ we need to encode a {}``function basis''
$\left\{ \tilde{f}_{b_{j}^{i}}\at_{e^{i}}\right\} $, where $\left\{ b_{j}^{i}\right\} $
is a basis for $R^{G_{e^{i}}}$. Since for different $e^{i}$'s the
spaces $R^{G_{e^{i}}}$ need not even overlap, we now have only function-chunks,
indexed by different $R$-elements for each edge, and no function
on the whole of $\Gamma$ to extract boundary conditions from. Fortunately,
algebra is generous and this complication turns out to be solvable.

\medskip{}

\subsection{\label{sub:construction-proof}Method}

We now present the actual construction procedure. Assume we have a
representation $R$ of a group $G$ acting on the quantum graph $\Gamma=\left(E,V\right)$,
and we have chosen representatives $\left\{ \tilde{e}^{i}\right\} _{i=1}^{I}$
for the orbits $\nicefrac{E}{G}$, and likewise $\left\{ \tilde{v}_{k}\right\} _{k=1}^{K}$
for $\nicefrac{V}{G}$. We have also chosen an ordered basis $B=\left(b_{j}\right)_{j=1}^{d}$
for $R$, and for each $i\in\left\{ 1..I\right\} $ another ordered
basis for $R$, $B^{i}=\left(b_{j}^{i}\right)_{j=1}^{d}$, such that
$\left\{ b_{j}^{i}\right\} _{j=1}^{d_{i}}$ is a basis for $R^{G_{\tilde{e}^{i}}}$
and each $b_{j}^{i}$ with $j>d_{i}=\dim R^{G_{\tilde{e}^{i}}}$ lies
in a nontrivial component of ${\rm Res}_{G_{\tilde{e}^{i}}}^{G}R$.

The quotient graph $\nicefrac{\Gamma}{R}$ obtained from these choices
is defined to have $\left\{ v_{k}\right\} _{k=1}^{K}$ as its set
of vertices, and $\left\{ e_{j}^{i}\right\} _{j=1..d_{i}}^{i=1..I}$
for edges, where each $e_{j}^{i}$ is of length $l_{\tilde{e}^{i}}$.
If $\tilde{e}^{i}$ connects $g\tilde{v}_{k}$ to $g'\tilde{v}_{k'}$
in $\Gamma$, then, for all $j$, $e_{j}^{i}$ connects $v_{k}$ to
$v_{k'}$ in $\nicefrac{\Gamma}{R}$. We shall assume, by adding {}``dummy''
vertices if needed, that $G$ does not carry any vertex in $V$ to
one of its neighbors. This serves three purposes:
\begin{enumerate}
\item It means that $\nicefrac{\Gamma}{R}$ has no loops; i.e., that $k\neq k'$
in the notation above. This allows us to speak of $f\at_{e_{j}^{i}}\left(v\right)$,
the value of $e_{j}^{i}$ at $v$, without confusion regarding which
end of $e_{j}^{i}$ is meant.
\item It assures that an edge is not transformed onto itself in the opposite
direction, in which case we would have had to take only half of the
edge as a representative for its orbit.
\item It assures that the fixed points of each $g\in G$ are either entire
edges, or vertices.
\end{enumerate}
Note that in order that $G$ still act on the graph, the dummy vertices
are to be added in accordance with its action, i.e., if a vertex is
placed at $x\in\left(0,l_{\tilde{e}}\right)$ along $\tilde{e}$,
one should also be placed at $x$ along $g\tilde{e}$, for every $g\in G$.

We can now define $\Psi$ on $\mathrm{Hom}_{\mathbb{C}G}\left(R,\mathcal{H}\left(\Gamma\right)\right)$:\[
\left(\Psi\tilde{f}\right)\at_{e_{j}^{i}}\overset{{\scriptstyle {\scriptscriptstyle def}}}{\equiv}\tilde{f}_{b_{j}^{i}}\at_{\tilde{e}^{i}}\:,\]
and it is clear that $\Psi$ does intertwine the Laplacians. We would
like to determine vertex conditions on $\nicefrac{\Gamma}{R}$ that
will ensure that $\Psi$ is into and onto $\mathcal{H}\left(\nicefrac{\Gamma}{R}\right)$.
This will require nothing more than linear algebra, and we start by
rephrasing \eqref{eq:edge-to-vector} basis-wise. We make the following
convention: an expression in bold is to be understood as a row vector
of length $d$, where the $\#$-symbol indicates the place of the
index; e.g., $\boldsymbol{\tilde{f}_{b_{\#}}\big|_{e}}$ stands for
$\left(\tilde{f}_{b_{1}}\big|_{e}\,,\,\ldots\,,\,\tilde{f}_{b_{d}}\big|_{e}\right)$.

Consider $\tilde{f}_{r}\at_{g\tilde{e}^{i}}$, an arbitrary function
in the family $\tilde{f}$ evaluated on an arbitrary edge. Write $r\in R$
as $\boldsymbol{b}\cdot\boldsymbol{\alpha}$, where $\boldsymbol{b}=\left(b_{1},\ldots,b_{d}\right)\in M_{1\times d}\left(R\right)$
and $\boldsymbol{\alpha}\in M_{d\times1}\left(\mathbb{C}\right)$.
$r=\boldsymbol{b}\cdot\boldsymbol{\alpha}$ implies $gr=\boldsymbol{b_{\#}^{i}}\cdot\left[\rho_{R}\left(g\right)\right]_{B^{i}}^{B}\cdot\boldsymbol{\alpha}$,
and therefore, by \eqref{eq:edge-to-vector} we have $\tilde{f}_{r}\at_{g\tilde{e}^{i}}=\tilde{f}_{\boldsymbol{b_{\#}^{i}}\left[\rho_{R}\left(g^{-1}\right)\right]_{B^{i}}^{B}\boldsymbol{\alpha}}\at_{\tilde{e}^{i}}$.
Linearity now implies \[
\tilde{f}_{r}\at_{g\tilde{e}^{i}}\equiv\tilde{f}_{\boldsymbol{b_{\#}^{i}}\left[g^{-1}\right]_{B^{i}}^{B}\boldsymbol{\alpha}}\at_{\tilde{e}^{i}}\equiv\boldsymbol{\tilde{f}_{b_{\#}^{i}}\at_{\tilde{e}^{i}}}\cdot\left[g^{-1}\right]_{B^{i}}^{B}\cdot\boldsymbol{\alpha}\equiv\boldsymbol{\left(\Psi\tilde{f}\right)\at_{e_{\#}^{i}}}\cdot\left[g^{-1}\right]_{B^{i}}^{B}\cdot\boldsymbol{\alpha}\,\,,\]
 where $\rho_{R}$ is understood, $\boldsymbol{\tilde{f}_{b_{\#}^{i}}\at_{\tilde{e}^{i}}}=\left(\tilde{f}_{b_{1}^{i}}\at_{\tilde{e}^{i}},\ldots,\tilde{f}_{b_{d}^{i}}\at_{\tilde{e}^{i}}\right)$,
and \[
\boldsymbol{\left(\Psi\tilde{f}\right)\at_{e_{\#}^{i}}}=\left(\left(\Psi\tilde{f}\right)\at_{e_{1}^{i}},\ldots,\left(\Psi\tilde{f}\right)\at_{e_{d_{i}}^{i}},0,\ldots,0\right)\,,\]
since for $j>d_{i}$ we have seen that $\tilde{f}_{b_{j}^{i}}\at_{\tilde{e}^{i}}\equiv0$,
and we therefore did not include the corresponding $e_{j}^{i}$ edge
in $\nicefrac{\Gamma}{R}$ (it is {}``dead'' - see footnote \ref{fn:dead-edges}).
We now see that for $f\in\mathcal{H}\left(\nicefrac{\Gamma}{R}\right)$
the inverse of $\Psi$ must be given by:\[
\left(\Psi^{-1}f\right)_{\mathbf{b}\cdot\boldsymbol{\alpha}}\at_{g\tilde{e}^{i}}\equiv\boldsymbol{f\at_{e_{\#}^{i}}}\cdot\left[g^{-1}\right]_{B^{i}}^{B}\cdot\boldsymbol{\alpha}\,,\]
(again $\left[g^{-1}\right]_{B^{i}}^{B}$ stands for $\left[\rho_{R}\left(g^{-1}\right)\right]_{B^{i}}^{B}$),
so we need to establish that the r.h.s is independent in the choice
of $g$. We observe that if $g$ and $g'$ are two possible choices
then $g^{-1}g'\in G_{\tilde{e}^{i}}$, and by the construction of
$B^{i}$ we have $\left[g^{-1}g'\right]_{B^{i}}=\left(\begin{array}{c|c}
I_{d_{i}} & 0\\
\hline 0 & *\end{array}\right)$. As we have agreed that $\boldsymbol{f\at_{e_{\#}^{i}}}=\left(f\at_{e_{1}^{i}},\ldots,f\at_{e_{d_{i}}^{i}},0,\ldots,0\right)$,
we have $\boldsymbol{f\at_{e_{\#}^{i}}}\cdot\left[g^{-1}g'\right]_{B^{i}}=\boldsymbol{f\at_{e_{\#}^{i}}}$
and thus \[
\boldsymbol{f\at_{e_{\#}^{i}}}\cdot\left[g'^{-1}\right]_{B^{i}}^{B}=\boldsymbol{f\at_{e_{\#}^{i}}}\cdot\left[g^{-1}g'\right]_{B^{i}}\cdot\left[g'^{-1}\right]_{B^{i}}^{B}=\boldsymbol{f\at_{e_{\#}^{i}}}\cdot\left[g^{-1}\right]_{B^{i}}^{B}\:,\]
establishing that $\Psi^{-1}$ is well defined.

\medskip{}

We can now determine matrices $A_{v_{k}}$ and $B_{v_{k}}$ for the
vertex $v_{k}$ from the matrices $A_{\tilde{v}_{k}}$, $B_{\tilde{v}_{k}}$
of the vertex $\tilde{v}_{k}$. Assume that the edges entering $\tilde{v}_{k}$
are $g_{1}\tilde{e}^{\nu_{1}},\ldots,g_{n}\tilde{e}^{\nu_{n}}$ (where
$n=d_{\tilde{v}_{k}}$), so that a function $f$ on $\Gamma$ satisfies
the vertex conditions at $\tilde{v}_{k}$ when \[
A_{\tilde{v}_{k}}\cdot f\at_{\tilde{v}_{k}}+B_{\tilde{v}_{k}}\cdot f'\at_{\tilde{v}_{k}}=0\:,\]
 where we recall from section \ref{sec:graphs} that\[
\begin{array}{cclcrc}
f\at_{\tilde{v}_{k}} & = & \Big(\, f\at_{g_{1}\tilde{e}^{\nu_{1}}}\left(\tilde{v}_{k}\right) & \ldots & f\at_{g_{n}\tilde{e}^{\nu_{n}}}\left(\tilde{v}_{k}\right)\,\Big)^{T}\\
f'\at_{\tilde{v}_{k}} & = & \Big(\, f'\at_{g_{1}\tilde{e}^{\nu_{1}}}\left(\tilde{v}_{k}\right) & \ldots & f'\at_{g_{n}\tilde{e}^{\nu_{n}}}\left(\tilde{v}_{k}\right)\,\Big)^{T} & .\end{array}\]

$\tilde{f}\in\mathrm{Hom}_{\mathbb{C}G}\left(R,\mathcal{H}\left(\Gamma\right)\right)$
means that $\tilde{f}_{r}$ satisfies the conditions at $\tilde{v}_{k}$
for all $r\in R$, which happens iff the basis functions $\left\{ \tilde{f}_{b_{j}}\right\} _{j=1}^{d}$
satisfy them. Thus, if we define the $n\times d$ matrix\begin{eqnarray*}
\boldsymbol{\tilde{f}_{b}\at_{\tilde{v}_{k}}} & = & \left(\begin{smallmatrix}\tilde{f}_{b_{1}}\at_{g_{1}\tilde{e}^{\nu_{1}}}\left(\tilde{v}_{k}\right) & \cdots & \tilde{f}_{b_{d}}\at_{g_{1}\tilde{e}^{\nu_{1}}}\left(\tilde{v}_{k}\right)\\
\vdots & \ddots & \vdots\\
\tilde{f}_{b_{1}}\at_{g_{n}\tilde{e}^{\nu_{n}}}\left(\tilde{v}_{k}\right) & \cdots & \tilde{f}_{b_{d}}\at_{g_{n}\tilde{e}^{\nu_{n}}}\left(\tilde{v}_{k}\right)\end{smallmatrix}\right)\\
 & = & \left(\begin{matrix}\tilde{f}_{b_{1}}\at_{\tilde{v}_{k}} & \cdots & \tilde{f}_{b_{d}}\at_{\tilde{v}_{k}}\end{matrix}\right)=\left(\begin{smallmatrix}\boldsymbol{\tilde{f}_{b_{\#}}\at_{g_{1}\tilde{e}^{\nu_{1}}}\left(\tilde{v}_{k}\right)}\\
\vdots\\
\boldsymbol{\tilde{f}_{b_{\#}}\at_{g_{n}\tilde{e}^{\nu_{n}}}\left(\tilde{v}_{k}\right)}\end{smallmatrix}\right)\:,\end{eqnarray*}
 and analogously $\boldsymbol{\tilde{f}_{b}'\at_{\tilde{v}_{k}}}$,
then we need only check that \begin{equation}
A_{\tilde{v}_{k}}\cdot\boldsymbol{\tilde{f}_{b}\at_{\tilde{v}_{k}}}+B_{\tilde{v}_{k}}\cdot\boldsymbol{\tilde{f}_{b}'\at_{\tilde{v}_{k}}}=0_{n\times d}\:.\label{eq:check-vertex-on-basis}\end{equation}
In addition, we note that if the boundary conditions are met by $\tilde{f}$
at $\tilde{v}_{k}$, then they are also met at any vertex in the orbit
$G\cdot\tilde{v}_{k}$, since $G$ is assumed to preserve boundary
conditions.

For a $n\times m$ matrix $X=\left(\left(x_{ij}\right)\right)$ we
define its row-wise-vectorization to be the $nm\times1$ matrix \[
\rv\, X\overset{{\scriptscriptstyle def}}{=}\left(\begin{array}{c}
\left(x_{11},\ldots,x_{1m}\right)^{T}\\
\vdots\\
\left(x_{n1},\ldots,x_{nm}\right)^{T}\end{array}\right)=\left(x_{11},x_{12},\ldots,x_{1m},x_{21},\ldots\ldots,x_{nm}\right)^{T}\:.\]
 Vectorization behaves quite nicely under multiplication. Specifically,
$\rv\left(A\cdot B\cdot C\right)=\left(A\otimes C^{T}\right)\cdot\rv\, B$,
which allows us to write \eqref{eq:check-vertex-on-basis} as\begin{equation}
\left(A_{\tilde{v}_{k}}\otimes I_{d}\right)\cdot\rv\,\boldsymbol{\tilde{f}_{b}\at_{\tilde{v}_{k}}}+\left(B_{\tilde{v}_{k}}\otimes I_{d}\right)\cdot\rv\,\boldsymbol{\tilde{f}_{b}'\at_{\tilde{v}_{k}}}=0_{nd\times1}\,\,.\label{eq:check-basis-vectorized}\end{equation}
Recalling that $\boldsymbol{\tilde{f}_{b_{\#}}\at_{g_{i}\tilde{e}^{\nu_{i}}}}=\boldsymbol{\Psi\tilde{f}\at_{e_{\#}^{\nu_{i}}}}\cdot\left[g_{i}^{-1}\right]_{B^{\nu_{i}}}^{B}$,
we have\begin{gather*}
\rv\,\boldsymbol{\tilde{f}_{b}\at_{\tilde{v}_{k}}}=\left(\begin{array}{c}
\boldsymbol{\tilde{f}_{b_{\#}}\at_{g_{1}\tilde{e}^{\nu_{1}}}\left(\tilde{v}_{k}\right)}^{T}\\
\vdots\\
\boldsymbol{\tilde{f}_{b_{\#}}\at_{g_{n}\tilde{e}^{\nu_{n}}}\left(\tilde{v}_{k}\right)}^{T}\end{array}\right)=\left(\begin{array}{c}
\left(\left[g_{1}^{-1}\right]_{B^{\nu_{1}}}^{B}\right)^{T}\cdot\left(\boldsymbol{\Psi\tilde{f}\at_{e_{\#}^{\nu_{1}}}\left(v_{k}\right)}\right)^{T}\\
\vdots\\
\left(\left[g_{n}^{-1}\right]_{B^{\nu_{n}}}^{B}\right)^{T}\cdot\left(\boldsymbol{\Psi\tilde{f}\at_{e_{\#}^{\nu_{n}}}\left(v_{k}\right)}\right)^{T}\end{array}\right)\\
={\rm diag}\left(\left[g_{1}^{-1}\right]_{B^{\nu_{1}}}^{B},\ldots,\left[g_{n}^{-1}\right]_{B^{\nu_{n}}}^{B}\right)^{T}\cdot\rv\left(\begin{array}{c}
\boldsymbol{\Psi\tilde{f}\at_{e_{\#}^{\nu_{1}}}\left(v_{k}\right)}\\
\vdots\\
\boldsymbol{\Psi\tilde{f}\at_{e_{\#}^{\nu_{n}}}\left(v_{k}\right)}\end{array}\right)\end{gather*}
and likewise for $\rv\,\boldsymbol{\tilde{f}_{b}'\at_{\tilde{v}_{k}}}$.
But now, the last vector is almost $\Psi\tilde{f}\at_{v_{k}}$, the
vector of values of $\left(\Psi\tilde{f}\right)$ at $v_{k}$! Only
two changes need to be made: first, if the edges entering $v_{k}$
are $e_{1}^{\mu_{1}},\ldots,e_{d_{\mu_{1}}}^{\mu_{1}},e_{1}^{\mu_{2}},\ldots,e_{d_{\mu_{m}}}^{\mu_{m}}$,
then by definition $\left\{ \mu_{i}\right\} _{i=1}^{m}=\left\{ \nu_{i}\right\} _{i=1}^{n}$
as sets; however, the $\mu_{i}$ are distinct, whereas in general,
repetitions can occur among the $\nu_{i}$ (i.e., two edges in $E_{\tilde{v}_{k}}$
might belong to the same $G$-orbit). Second, as in all our expressions
there might be {}``dead'' edges, $e_{j}^{\mu_{i}}$ with $j>d_{\mu_{i}}$,
which do not really appear in the quotient graph (note, however, that
neither of the problems can occur when the action of $G$ is free).
We shall deal with these two inconveniences at once: we define the
$n\times m$ matrix $\left(\Theta'\right)_{ij}=\begin{cases}
1 & \nu_{i}=\mu_{j}\\
0 & otherwise\end{cases}$, and then take $\Theta$ to be the $nd\times d_{v_{k}}$ matrix obtained
by removing from $\left(\Theta'\otimes I_{d}\right)$ the columns
$\left\{ \left(i-1\right)\cdot d+j\right\} _{{1\leq i\leq m\atop d_{\mu_{i}}<j\leq d}}$;
these are the columns which would have been multiplied by a {}``dead''
edge in $\left(\begin{smallmatrix}\Psi\tilde{f}\big|_{e_{1}^{\mu_{1}}} & \ldots & \Psi\tilde{f}\big|_{e_{d}^{\mu_{1}}} & \Psi\tilde{f}\big|_{e_{1}^{\mu_{2}}} & \ldots & \ldots & \Psi\tilde{f}\big|_{e_{d}^{\mu_{m}}}\end{smallmatrix}\right)^{T}$
. We now have \[
\rv\left(\begin{matrix}\boldsymbol{\Psi\tilde{f}\at_{e_{\#}^{\nu_{1}}}\left(v_{k}\right)}\\
\vdots\\
\boldsymbol{\Psi\tilde{f}\at_{e_{\#}^{\nu_{n}}}\left(v_{k}\right)}\end{matrix}\right)=\Theta\cdot\left(\begin{smallmatrix}\Psi\tilde{f}\big|_{e_{1}^{\mu_{1}}}\left(v_{k}\right) & \ldots & \Psi\tilde{f}\big|_{e_{d_{\mu_{1}}}^{\mu_{1}}}\left(v_{k}\right) & \ldots\ldots & \Psi\tilde{f}\big|_{e_{d_{\mu_{m}}}^{\mu_{m}}}\left(v_{k}\right)\end{smallmatrix}\right)^{T}=\Theta\cdot\Psi\tilde{f}\at_{v_{k}}\:,\]
 and we can thus define\begin{eqnarray}
A_{v_{k}} & = & \left(A_{\tilde{v}_{k}}\otimes I_{d}\right)\cdot\mathfrak{G}\cdot\Theta\label{eq:final-boundary-conditions-1}\\
B_{v_{k}} & = & \left(B_{\tilde{v}_{k}}\otimes I_{d}\right)\cdot\mathfrak{G}\cdot\Theta\label{eq:final-boundary-conditions-2}\end{eqnarray}
where $\mathfrak{G}={\rm diag}\left(\left[g_{1}^{-1}\right]_{B^{\nu_{1}}}^{B},\ldots,\left[g_{n}^{-1}\right]_{B^{\nu_{n}}}^{B}\right)^{T}$,
and finally rewrite \eqref{eq:check-basis-vectorized} as\[
A_{v_{k}}\cdot\Psi\tilde{f}\at_{v_{k}}+B_{v_{k}}\cdot\Psi\tilde{f}'\at_{v_{k}}=0\:.\]
These vertex conditions on $\Psi\tilde{f}$ at $v_{k}$ are equivalent
to $\tilde{f}_{r}$ satisfying the vertex conditions at $\tilde{v}_{k}$
for all $r\in R$, and therefore also on the entire orbit $G\cdot\tilde{v}_{k}$.
If we repeat this process for each $k=1\ldots K$, we indeed obtain
boundary conditions on $\nicefrac{\Gamma}{R}$ which are satisfied
by $\Psi\tilde{f}$ exactly when $\tilde{f}\in\mathrm{Hom}_{\mathbb{C}G}\left(R,\mathcal{H}\left(\Gamma\right)\right)$.

If the action of $G$ is free, then $\Theta$ is just a permutation
matrix (we can even order $E_{v_{k}}$ so that $\Theta=I$), but in
the general case $\Theta$ might be non-square (explicitly, it is
of size $nd\times d_{v_{k}}$, where $d_{v_{k}}=\sum_{i=1}^{m}d_{\mu_{i}}\leq md\leq nd$).
When this occurs, the matrices $A_{v_{k}}$ and $B_{v_{k}}$ we have
obtained are not square matrices, and we therefore obtain a quotient
which is only a generalized quantum graph. Nevertheless, as the matrices
$A_{v_{k}}$ and $B_{v_{k}}$ serve only to represent the system of
equations $A_{v_{k}}\cdot f\at_{v_{k}}+B_{v_{k}}\cdot f'\at_{v_{k}}=0$,
we can perform elementary row operations on the $nd\times2d_{v_{k}}$
matrix $\left(A_{v_{k}}\,\middle|\, B_{v_{k}}\right)$ without changing
the boundary conditions at $v_{k}$, and thus perhaps reduce the number
of rows of $\left(A_{v_{k}}\,\middle|\, B_{v_{k}}\right)$. In the
case that $\mathrm{rank}\left(A_{v_{k}}\,\middle|\, B_{v_{k}}\right)\leq d_{v_{k}}$,
we can reduce the matrices $A_{v_{k}}$ and $B_{v_{k}}$ to squares
ones, and if this holds for all $k$ then we actually have a proper
quantum graph. If it further happens that $\mathrm{rank}\left(A_{v_{k}}\,\middle|\, B_{v_{k}}\right)=d_{v_{k}}$
for all $k$, then the quotient graph is also exact. We now show sufficient
conditions for this to happen.
\begin{prop}
\label{pro:rank-of-SA-quotient}If there exist $\omega\in\mathbb{C}^{\times}$
and $M\in\mathrm{GL}_{d_{\tilde{v}_{k}}}\left(\mathbb{C}\right)$
such that $\left(A_{\tilde{v}_{k}}\,\middle|\, B_{\tilde{v}_{k}}\right)$
is row-equivalent to $\left(\omega\left(M-I\right)\,\middle|\, M+I\right)$,
then $\mathrm{rank}\left(A_{v_{k}}\,\middle|\, B_{v_{k}}\right)=d_{v_{k}}$.\end{prop}
\begin{rem*}
We recall from section \ref{sec:graphs} that this condition holds
for all $k$ when $\Gamma$'s Laplacian is self-adjoint. Therefore,
in this case $\nicefrac{\Gamma}{R}$ is exact, as stated in theorem
\ref{thm:construction-thm}.\end{rem*}
\begin{proof}
Denote $\tilde{v}=\tilde{v}_{k}$, $v=v_{k}$, and recall that $E_{\tilde{v}}=\left\{ g_{i}\tilde{e}^{\nu_{i}}\right\} _{i=1}^{n}$
is the set of edges entering $\tilde{v}$. Assume, by reordering if
necessary, that $\nu_{i}=\mu_{i}$ for $1\leq i\leq m$, i.e., that
$\left\{ g_{i}\tilde{e}^{\nu_{i}}\right\} _{i=1}^{m}$ are representatives
for the $G_{\tilde{v}}$-orbits in $E_{\tilde{v}}$. Denote $\tilde{\varepsilon}^{i}=g_{i}\tilde{e}^{\nu_{i}}=g_{i}\tilde{e}^{\mu_{i}}$
(where $1\leq i\leq m$), and note that $G_{\tilde{\varepsilon}^{i}}$
is conjugate to $G_{\tilde{e}^{\mu_{i}}}$. The action of $G_{\tilde{v}}$
on $E_{\tilde{v}}$ gives rise to a representation $\mathbb{C}\left[E_{\tilde{v}}\right]$
of $G_{\tilde{v}}$, and the $G_{\tilde{v}}$-set isomorphism $E_{\tilde{v}}=\coprod_{i=1}^{m}G_{\tilde{v}}\cdot\tilde{\varepsilon}^{i}\cong\coprod_{i=1}^{m}{}^{G_{\tilde{v}}}\!/\!/\!_{G_{\tilde{\varepsilon}^{i}}}$
translates to an isomorphism of $G_{\tilde{v}}$-representations:
\[
\mathbb{C}\left[E_{\tilde{v}}\right]\cong\bigoplus_{i=1}^{m}\mathbb{C}\left[^{G_{\tilde{v}}}\!/\!/\!_{G_{\tilde{\varepsilon}^{i}}}\right]\cong\bigoplus_{i=1}^{m}\ind_{G_{\tilde{\varepsilon}^{i}}}^{G_{\tilde{v}}}\mathbf{1}_{G_{\tilde{\varepsilon}^{i}}}\:.\]
Here $\mathbf{1}_{G}$ denotes the trivial representation of a group
$G$, but we shall also use it to denote its character. We now see
that \begin{multline}
\left\langle \chi_{\mathbb{C}\left[E_{\tilde{v}}\right]},\chi_{R}\right\rangle _{G_{\tilde{v}}}=\left\langle \chi_{\bigoplus_{i=1}^{m}\ind_{G_{\tilde{\varepsilon}^{i}}}^{G_{\tilde{v}}}\mathbf{1}_{G_{\tilde{\varepsilon}^{i}}}},\chi_{R}\right\rangle _{G_{\tilde{v}}}=\sum_{i=1}^{m}\left\langle \ind_{G_{\tilde{\varepsilon}^{i}}}^{G_{\tilde{v}}}\mathbf{1}_{G_{\tilde{\varepsilon}^{i}}},\chi_{R}\right\rangle _{G_{\tilde{v}}}\\
=\sum_{i=1}^{m}\left\langle \mathbf{1}_{G_{\tilde{\varepsilon}^{i}}},\chi_{R}\right\rangle _{G_{\tilde{\varepsilon}^{i}}}=\sum_{i=1}^{m}\dim R^{G_{\tilde{\varepsilon}^{i}}}=\sum_{i=1}^{m}\dim R^{G_{\tilde{e}^{\mu_{i}}}}=\sum_{i=1}^{m}d_{\mu_{i}}=d_{v}\:.\label{eq:R1-significance}\end{multline}
We return to the matrices $\left(A_{v}\,\middle|\, B_{v}\right)\in M_{nd\times2d_{v}}\left(\mathbb{C}\right)$
and $\left(A_{\tilde{v}}\,\middle|\, B_{\tilde{v}}\right)\in M_{n\times2n}\left(\mathbb{C}\right)$.
For $f\in\mathcal{H}\left(\Gamma\right)$, the action of $G_{\tilde{v}}$
on $E_{\tilde{v}}$ induces a permutation action of $G_{\tilde{v}}$
on the entries of $f\at_{\tilde{v}}=\left(f\at_{\tilde{e}}\left(\tilde{v}\right)\right)_{\tilde{e}\in E_{\tilde{v}}}$,
and exactly the same action is induced on the entries of $f'\Big|_{\tilde{v}}$.
Thus, the space $\mathbb{C}^{2n}$ of possible values and derivatives
at $\tilde{v}$ has naturally the structure of the $G_{\tilde{v}}$-representation
$\mathbb{C}\left[E_{\tilde{v}}\right]\oplus\mathbb{C}\left[E_{\tilde{v}}\right]$.
Furthermore, as by assumption $G$ preserves the boundary conditions,
$\ker\left(A_{\tilde{v}}\,\middle|\, B_{\tilde{v}}\right)\subseteq\mathbb{C}^{2n}$
is a sub-$G_{\tilde{v}}$-representation of $\mathbb{C}^{2n}\cong\mathbb{C}\left[E_{\tilde{v}}\right]\oplus\mathbb{C}\left[E_{\tilde{v}}\right]$.
We observe that the encoding and decoding processes are {}``rigid'',
in the sense that for $x\in\left[0,l_{\tilde{e}^{i}}\right]$ it suffices
to know $\left\{ \tilde{f}_{r}\at_{\tilde{e}^{i}}\left(x\right)\right\} _{r\in R}$
to determine $\left\{ \Psi\tilde{f}\at_{e_{j}^{i}}\left(x\right)\right\} _{j=1..d_{i}}$,
and vice versa. Likewise, $\left\{ \tilde{f}_{r}\at_{\tilde{v}}\right\} _{r\in R}$
and $\Psi\tilde{f}\at_{v}$ determine one another, and the same goes
for the corresponding derivatives. This means that in the commutative
diagram\[ \xymatrix@R=9pt@C=2pt{**[l]\mathrm{Hom}_{\mathbb{C}G}\left(R,\mathcal{H}\left(\Gamma\right)\right)\ar[r]\ar[ddd]<-12pt>_{\Psi} & **[r]\mathrm{Hom}_{\mathbb{C}G_{\tilde{v}}}\left(R,\ker\left(A_{\tilde{v}}\,\middle|\, B_{\tilde{v}}\right)\right)\ar[ddd]<12pt>^{\psi}\\ {\scriptstyle \tilde{f}}\ar@{|->}[r]\ar@{|->}[d] & **[l]{\scriptstyle \left(r\mapsto\left(\tilde{f}_{r}\big|_{\tilde{v}},\tilde{f}'_{r}\big|_{\tilde{v}}\right)\right)}\ar@{|->}[d]<-26pt>\\ {\scriptstyle \Psi\tilde{f}}\ar@{|->}[r] & **[l]{\scriptstyle \left(\Psi\tilde{f}\big|_{v},\left(\Psi\tilde{f}\right)'\big|_{v}\right)}\\ **[l]\mathcal{H}\left(\nicefrac{\Gamma}{R}\right)\ar[r] & **[r]\ker\left(A_{v}\,\middle|\, B_{v}\right)} \]the
map $\psi$, which is this {}``local'' encoding, is in fact an isomorphism.
This gives us \begin{equation}
\mathrm{null}\left(A_{v}\,\middle|\, B_{v}\right)=\left\langle \chi_{R},\chi_{\ker\left(A_{\tilde{v}}\,\middle|\, B_{\tilde{v}}\right)}\right\rangle _{G_{\tilde{v}}}\:,\label{eq:R2-significance}\end{equation}
so that by \eqref{eq:R1-significance}\[
\mathrm{rank}\left(A_{v}\,\middle|\, B_{v}\right)=\left\langle 2\chi_{\mathbb{C}\left[E_{\tilde{v}}\right]}-\chi_{\ker\left(A_{\tilde{v}}\,\middle|\, B_{\tilde{v}}\right)},\chi_{R}\right\rangle _{G_{\tilde{v}}}\:.\]
We therefore have\begin{eqnarray*}
\qquad\quad\mathrm{rank}\left(A_{v}\,\middle|\, B_{v}\right)=d_{v} & \qquad\Leftrightarrow\quad & \left\langle \chi_{\mathbb{C}\left[E_{\tilde{v}}\right]}-\chi_{\ker\left(A_{\tilde{v}}\,\middle|\, B_{\tilde{v}}\right)},\chi_{R}\right\rangle _{G_{\tilde{v}}}=0\:,\end{eqnarray*}
and the last equality holds for all representations $R$ of $G$ if
and only if $\ind_{G_{\tilde{v}}}^{G}\mathbb{C}\left[E_{\tilde{v}}\right]\cong\ind_{G_{\tilde{v}}}^{G}\ker\left(A_{\tilde{v}}\,\middle|\, B_{\tilde{v}}\right)$.
In particular, this happens if $\mathbb{C}\left[E_{\tilde{v}}\right]$
and $\ker\left(A_{\tilde{v}}\,\middle|\, B_{\tilde{v}}\right)$ are
isomorphic $G_{\tilde{v}}$-representations, which we now show to
follow from our assumptions.

Observe that $\xi:\mathbb{C}\left[E_{\tilde{v}}\right]\oplus\mathbb{C}\left[E_{\tilde{v}}\right]\rightarrow\mathbb{C}\left[E_{\tilde{v}}\right]$,
defined by $\xi\left(\underline{a},\underline{b}\right)=\omega\underline{a}-\underline{b}$
is a homomorphism of $G_{\tilde{v}}$-representations, and recall
that $\ker\left(A_{\tilde{v}}\,\middle|\, B_{\tilde{v}}\right)$ is
naturally embedded in $\mathbb{C}\left[E_{\tilde{v}}\right]\oplus\mathbb{C}\left[E_{\tilde{v}}\right]$.
When restricting $\xi$ to $\ker\left(A_{\tilde{v}}\,\middle|\, B_{\tilde{v}}\right)$
we obtain the desired isomorphism onto $\mathbb{C}\left[E_{\tilde{v}}\right]$,
since $\dim\ker\left(A_{\tilde{v}}\,\middle|\, B_{\tilde{v}}\right)=\mathrm{null}\left(\omega\left(M-I\right)\,\middle|\,\left(M+I\right)\right)=d_{\tilde{v}}=\dim\mathbb{C}\left[E_{\tilde{v}}\right]$,
and\[
\left(\underline{a},\underline{b}\right)\in\ker\left(\xi\Big|_{\ker\left(A_{\tilde{v}}\,\middle|\, B_{\tilde{v}}\right)}\right)\,\Rightarrow\,\left\{ \begin{array}{c}
\omega\left(M-I\right)\underline{a}+\left(M+I\right)\underline{b}=0\\
\omega\underline{a}-\underline{b}=0\end{array}\right\} \:\Rightarrow\:\left(\underline{a},\underline{b}\right)=0\:.\]

\end{proof}

\subsection{\label{sub:construction-remarks}Remarks}

\

\subsubsection{\label{sub:smooth-quotient}}

It seems of some interest to point out that the encoding process we
have described has actually nothing to do with eigenfunctions of the
Laplacian. The assumption that $\tilde{f}_{r}\in\mathcal{H}\left(\Gamma\right)$
was not used during the construction of the quotient, and as a result,
if no dummy vertices are introduced at the beginning of the construction,
then we actually have \begin{equation}
\Psi:\mathrm{Hom}_{\mathbb{C}G}\left(R,C^{^{\infty}}\negthickspace\left(\Gamma\right)\right)\overset{\cong}{\longrightarrow}C^{^{\infty}}\negthickspace\left(\nicefrac{\Gamma}{R}\right)\:.\label{eq:smooth-quotient}\end{equation}
If dummy vertices are added, and $\Gamma'$ is the graph obtained
from $\Gamma$ by their introduction, we obtain only $C^{^{\infty}}\negthickspace\left(\nicefrac{\Gamma}{R}\right)\cong\mathrm{Hom}_{\mathbb{C}G}\left(R,C^{^{\infty}}\negthickspace\left(\Gamma'\right)\right)$,
and unfortunately $C^{^{\infty}}\negthickspace\left(\Gamma'\right)\neq C^{^{\infty}}\negthickspace\left(\Gamma\right)$,
as was remarked in section \ref{sec:graphs}. We have introduced dummy
vertices in order to avoid loops and parallel edges, and also to ensure
that a subset of the edges can be taken as a fundamental domain for
the graph. Of these causes, only the last is unavoidable; one can
still carry out the construction (encumbering somewhat the notations),
even with loops or parallel edges, both in the original graph and
in the quotient. The only case in which the construction fails altogether,
and a dummy point must be introduced, is when a group element inverts
the direction of an edge, so that a fundamental domain must include
only half of the edge. Thus, in order for it to be possible to construct
by our method a {}``smooth quotient'', in the sense of \eqref{eq:smooth-quotient},
$\mathrm{fix}_{g}\Gamma$ must be a subgraph of $\Gamma$ for every
$g\in G$. We shall return to these observations in section \ref{sub:squareandtriangle}.

\subsubsection{}

If $G$ acts on $\Gamma$ and $R$ is a representation of $H\leq G$,
we can consider the composition of isomorphisms \[
\mathcal{H}\left(\nicefrac{\Gamma}{R}\right)\overset{\Psi^{-1}}{\longrightarrow}\mathrm{Hom}_{\mathbb{C}H}\left(R,\mathcal{H}\left(\Gamma\right)\right)\overset{\mathcal{F}}{\rightarrow}\mathrm{Hom}_{\mathbb{C}G}\left(\ind_{H}^{G}R,\mathcal{H}\left(\Gamma\right)\right)\overset{\Psi'}{\longrightarrow}\mathcal{H}\left(\nicefrac{\Gamma}{\ind_{H}^{G}R}\right)\]
where $\Psi$ and $\Psi'$ are the isomorphisms defined during the
constructions of $\nicefrac{\Gamma}{R}$ and $\nicefrac{\Gamma}{\ind_{H}^{G}R}$,
respectively, and $\mathcal{F}$ is the Frobenius isomorphism%
\footnote{Taking the induction to be the scalar extension $\ind_{H}^{G}R=\mathbb{C}G\otimes_{\mathbb{C}H}R$,
$\mathcal{F}$ is defined for decomposable tensors by $\left(\mathcal{F}\tilde{f}\right)_{g\otimes r}=g\cdot\tilde{f}_{r}$
(where $\tilde{f}\in\mathrm{Hom}_{\mathbb{C}H}\left(R,\mathcal{H}\left(\Gamma\right)\right)$,
$g\in G$, $r\in R$), and extends linearly to all of the tensor product.%
}. We obtain what is known as a \emph{transplantation} (see \cite{Buser,Berard})
between $\nicefrac{\Gamma}{R}$ and $\nicefrac{\Gamma}{\ind_{H}^{G}R}$,
an operator which construct functions on one graph as linear combinations
of segments of functions on the second graph. This is developed in
more details in \cite{BPbS}.

\subsubsection{}

It is natural to ask, for a graph $\Gamma$ whose Laplacian is self-adjoint,
whether the Laplacian on $\nicefrac{\Gamma}{R}$ is self-adjoint.
This turns out to depend on both the action of $G$ and the choices
of bases in the construction, and it is addressed for some cases in
\cite{BPbS}.

\subsubsection{}

Another natural question is the following: for a quantum graph $\Gamma$
acted upon by $G$, when does an irreducible representation $S$ of
$G$ appear in $\mathcal{H}\left(\Gamma\right)$?%
\footnote{This question, in the context of compact Lie groups acting on Riemannian
manifolds, is addressed in \cite{Brun}.%
} It is known that every quantum graph with edges whose Laplacian is
self-adjoint has a nonempty spectrum (see for example \cite{BolEnd}).
Therefore, if $\nicefrac{\Gamma}{S}$'s Laplacian is self-adjoint
then $S$ appears in $\mathcal{H}\left(\Gamma\right)$ iff $\nicefrac{\Gamma}{S}$
has edges, and by the construction method this happens iff for at
least one edge $e$ in $\Gamma$ the representation ${\rm Res}_{G_{e}}^{G}S$
has a nonempty trivial component, i.e., $\left\langle \chi_{S},\mathbf{1}\right\rangle _{G_{e}}\neq0$.
In particular, if $\Gamma$'s Laplacian is self adjoint, and $G$
acts freely on $\Gamma$, then a self-adjoint quotient can always
be obtained \cite{BPbS}, and each stabilizer has only the trivial
irreducible representation. Thus, every irreducible representation
of $G$ appears in $\mathcal{H}\left(\Gamma\right)$.

\section{Examples of isospectral quantum graphs}

\label{sec:isospectral_qg}We now demonstrate several applications
of the theory presented above which yield isospectral graphs. All
the examples below are direct consequences of the theorem or the corollary
presented in section \ref{sec:algebra}.

\begin{center}
\begin{figure}[!h]
\begin{centering}
\includegraphics[scale=0.4]{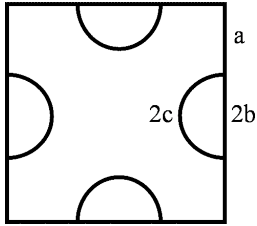} 
\par\end{centering}

\caption{A graph that obeys the dihedral symmetry of the square. The lengths
of some edges are marked.}

\label{fig:full_square} 
\end{figure}

\par\end{center}

Let $\Gamma$ be the graph given in figure \ref{fig:full_square}.
The lengths of the edges are determined by the parameters $a,b,c$
and it has Neumann boundary conditions at all vertices. $G=D_{4}$,
the dihedral group of the square, is a symmetry group of $\Gamma$.
Denote by $\tau$ the reflection of $\Gamma$ along the horizontal
axis and by $\sigma$ the rotation of $\Gamma$ counterclockwise by
$\nicefrac{\pi}{2}$. Then we can describe $G$ and some of its subgroups
$H_{1},\, H_{2},\, H_{3}\leq G$ by: \begin{eqnarray*}
G & = & \{e,\,\sigma,\,\sigma^{2},\,\sigma^{3},\,\tau,\,\tau\sigma,\,\tau\sigma^{2},\,\tau\sigma^{3}\}\\
H_{1} & = & \{e,\,\tau,\,\tau\sigma^{2},\,\sigma^{2}\}\\
H_{2} & = & \{e,\,\tau\sigma,\,\tau\sigma^{3},\,\sigma^{2}\}\\
H_{3} & = & \{e,\,\sigma,\,\sigma^{2},\,\sigma^{3}\}\end{eqnarray*}
Consider the following one dimensional representations of $H_{1}$,
$H_{2}$ and $H_{3}$ respectively: \begin{eqnarray}
R_{1} & : & \left\{ \begin{array}{llll}
e\mapsto\left(1\right), & \tau\mapsto\left(-1\right), & \tau\sigma^{2}\mapsto\left(1\right), & \sigma^{2}\mapsto\left(-1\right)\end{array}\right\} \label{eq:r1_rep}\\
R_{2} & : & \left\{ \begin{array}{llll}
e\mapsto\left(1\right), & \tau\sigma\mapsto\left(1\right), & \tau\sigma^{3}\mapsto\left(-1\right), & \sigma^{2}\mapsto\left(-1\right)\end{array}\right\} \label{eq:r2_rep}\\
R_{3} & : & \left\{ \begin{array}{llll}
e\mapsto\left(1\right), & \sigma\mapsto\left(i\right), & \sigma^{2}\mapsto\left(-1\right), & \sigma^{3}\mapsto\left(-i\right)\end{array}\right\} \label{eq:r3_rep}\end{eqnarray}
These representations fulfill the condition in corollary \ref{cor:sunadapair}:
$\ind_{H_{1}}^{G}R_{1}\cong\ind_{H_{2}}^{G}R_{2}\cong\ind_{H_{3}}^{G}R_{3}$
and thus we obtain that $\nicefrac{\Gamma}{R_{1}},\,\nicefrac{\Gamma}{R_{2}}$
and $\nicefrac{\Gamma}{R_{3}}$ are isospectral (figure \ref{fig:dihedral_triple}).

\begin{center}
\begin{figure}[!h]
\begin{centering}
\begin{minipage}[c]{0.3\columnwidth}%
\begin{center}
\includegraphics[scale=0.35]{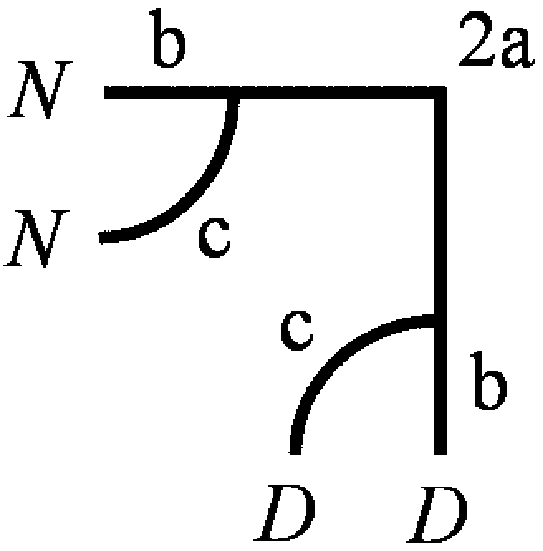} 
\par\end{center}

\begin{center}
(a) 
\par\end{center}%
\end{minipage}\hfill{}%
\begin{minipage}[c]{0.3\columnwidth}%
\begin{center}
\includegraphics[scale=0.35]{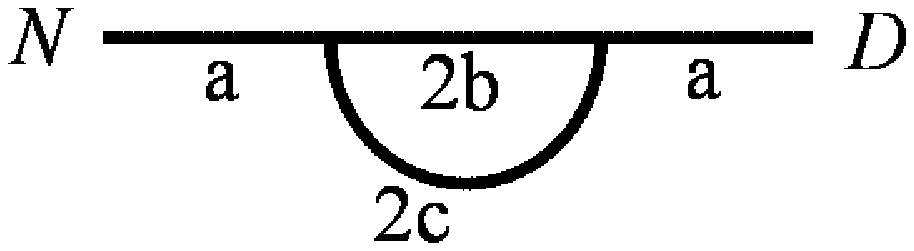} 
\par\end{center}

\medskip{}

\begin{center}
(b) 
\par\end{center}%
\end{minipage}\hfill{}\hfill{}\hfill{}%
\begin{minipage}[c][1\totalheight][s]{0.33\columnwidth}%
\begin{center}
\begin{minipage}[c]{0.4\columnwidth}%
\begin{center}
\includegraphics[clip,width=0.75\textwidth]{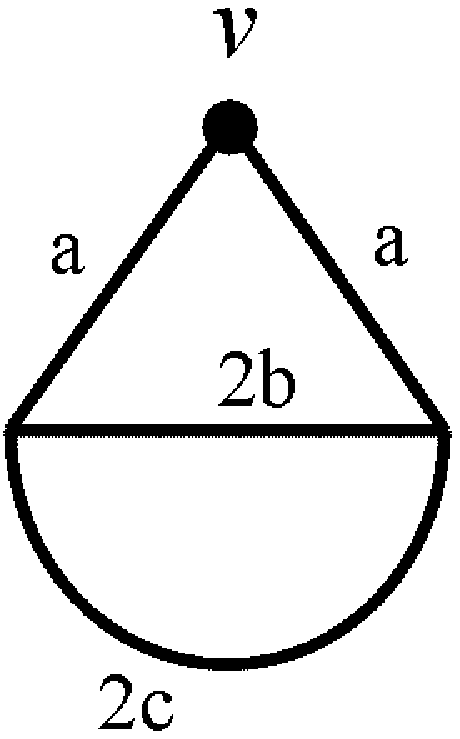} 
\par\end{center}%
\end{minipage}%
\begin{minipage}[c]{0.5\columnwidth}%
\begin{center}
$\begin{array}{l}
A_{v}=\left(\begin{smallmatrix}1 & i\\
0 & 0\end{smallmatrix}\right)\\
\\B_{v}=\left(\begin{smallmatrix}0 & 0\\
1 & -i\end{smallmatrix}\right)\end{array}$ 
\par\end{center}%
\end{minipage}
\par\end{center}

\begin{center}
(c) 
\par\end{center}%
\end{minipage}
\par\end{centering}

\caption{The three isospectral graphs $\nicefrac{\Gamma}{R_{1}},\nicefrac{\Gamma}{R_{2}},\nicefrac{\Gamma}{R_{3}}$.
Neumann boundary conditions are assumed if nothing else is specified.
$D$ stands for Dirichlet boundary conditions and $N$ for Neumann. }

\label{fig:dihedral_triple} 
\end{figure}

\par\end{center}

We now explain the process of building the graph $\nicefrac{\Gamma}{R_{1}}$.
First we give an intuition which suffices to obtain the quotient in
this case, and afterwards we strictly implement the method that is
described in section \ref{sub:construction-proof}. Going back to
(\ref{eq:quotient-definition}), we observe that the r.h.s. of it
is \begin{equation}
\mathrm{Hom}_{\mathbb{C}H_{1}}\left(R_{1},\mathcal{H}\left(\Gamma\right)\right)\cong\left(\mathcal{H}\left(\Gamma\right)\right)^{R_{1}}\,,\label{eq:hom-simple-isotip}\end{equation}
 where $\left(\mathcal{H}\left(\Gamma\right)\right)^{R_{1}}$ is the
$R_{1}$-isotypic component of $\mathcal{H}\left(\Gamma\right)$ (considered
as a $\mathbb{C}H_{1}$-module); the isomorphism is due to the fact
that $R_{1}$ is one-dimensional, hence irreducible. Let us study
the properties of $\tilde{f}\in\left(\mathcal{H}\left(\Gamma\right)\right)^{R_{1}}$.
We know (see \eqref{eq:r1_rep}) that $\tau\tilde{f}=-\tilde{f}$,
which means that $\tilde{f}$ is an anti-symmetric function with respect
to the horizontal reflection. We deduce that $\tilde{f}$ vanishes
on the fixed points of $\tau$ (marked with diamonds in figure \ref{fig:intuitive_quotient1}(a)).

\begin{center}
\begin{figure}[H]
\hfill{}%
\begin{minipage}[c][1\totalheight][t]{0.4\columnwidth}%
\begin{center}
\includegraphics[width=0.6\textwidth]{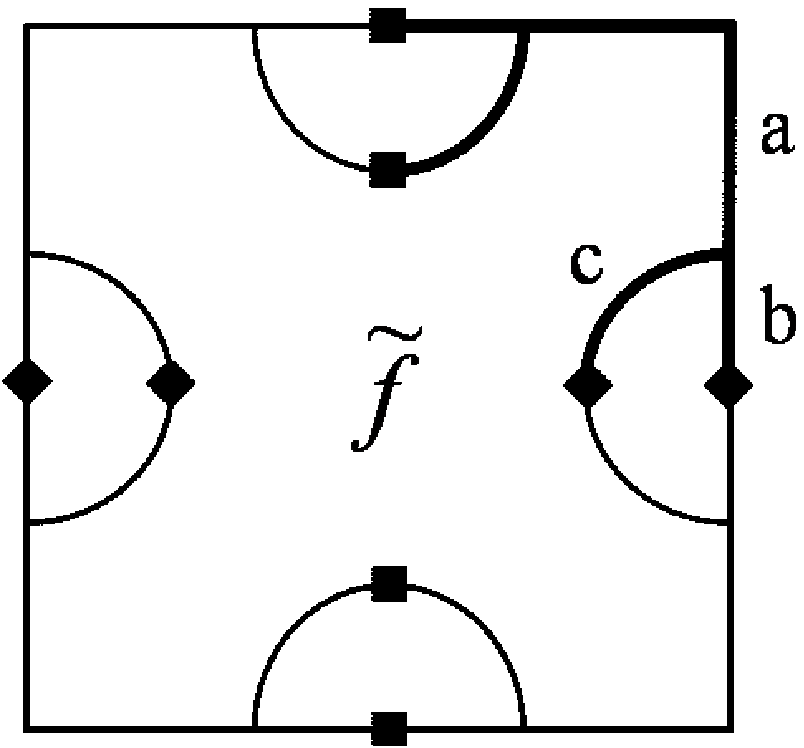} 
\par\end{center}

\begin{center}
(a) 
\par\end{center}%
\end{minipage}\hfill{}%
\begin{minipage}[c][1\totalheight][t]{0.35\columnwidth}%
\begin{center}
\includegraphics[width=0.6\textwidth]{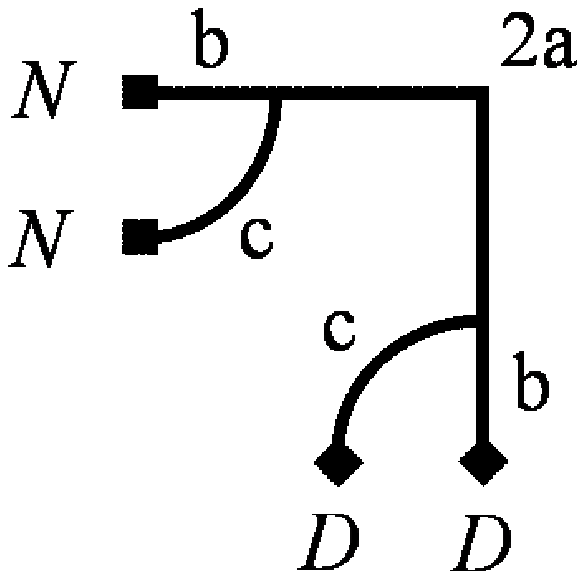} 
\par\end{center}

\begin{center}
(b) 
\par\end{center}%
\end{minipage}\hfill{}

\caption{(a) The information we have on $\tilde{f}\in\left(\mathcal{H}\left(\Gamma\right)\right)^{R_{1}}$.
Diamonds mark the vertices on which the function vanishes and squares
the vertices with zero derivative. (b) The quotient graph $\nicefrac{\Gamma}{R_{1}}$
which encodes this information. $D$ stands for Dirichlet boundary
conditions and $N$ for Neumann.}

\label{fig:intuitive_quotient1} 
\end{figure}

\par\end{center}

In a similar manner, we see that $\tilde{f}$ is symmetric with respect
to the vertical reflection since $\tau\sigma^{2}\tilde{f}=\tilde{f}$,
and therefore the derivative of $\tilde{f}$ must vanish at the corresponding
points (the squares in figure \ref{fig:intuitive_quotient1}(a)).
Furthermore, it is enough to know the values of $\tilde{f}$ restricted
to the first quadrant (the bold subgraph in figure \ref{fig:intuitive_quotient1}(a))
in order to deduce $\tilde{f}$ on the whole graph, using the known
action of the reflections, which follows from $\tilde{f}\in\left(\mathcal{H}\left(\Gamma\right)\right)^{R_{1}}$:
\begin{equation}
\tau\tilde{f}=-\tilde{f}\,,\,\,\tau\sigma^{2}\tilde{f}=\tilde{f}\,.\label{eq:R1_isotip_reflec}\end{equation}
 Our encoding is now complete and the quotient $\nicefrac{\Gamma}{R_{1}}$
is the subgraph which lies in the first quadrant, with the boundary
conditions of Dirichlet and Neumann in the appropriate locations as
was found for $\tilde{f}$ (figure \ref{fig:intuitive_quotient1}(b)).
The encoding is described by the map $\Psi:\mathrm{Hom}_{\mathbb{C}G}\left(R_{1},\mathcal{H}\left(\Gamma\right)\right)\rightarrow\mathcal{H}\left(\nicefrac{\Gamma}{R_{1}}\right)$
which is just the restriction map of functions in $\mathrm{Hom}_{\mathbb{C}G}\left(R_{1},\mathcal{H}\left(\Gamma\right)\right)\cong\mathcal{H}\left(\Gamma\right)^{R_{1}}$
to the mentioned subgraph. An important observation is that given
$f\in\mathcal{H}\left(\nicefrac{\Gamma}{R_{1}}\right)$ it is possible
to construct a unique function $\tilde{f}\in\left(\mathcal{H}\left(\Gamma\right)\right)^{R_{1}}$
(using \eqref{eq:R1_isotip_reflec}), whose restriction to the first
quadrant subgraph is $f$. It follows that $\Psi$ is invertible and
thus is an isomorphism. This ends the intuitive approach and we now
proceed to the rigorous derivation.

\begin{center}
\begin{figure}[!h]
\hfill{}%
\begin{minipage}[c][1\totalheight][t]{0.4\columnwidth}%
\begin{center}
\includegraphics[width=0.8\textwidth]{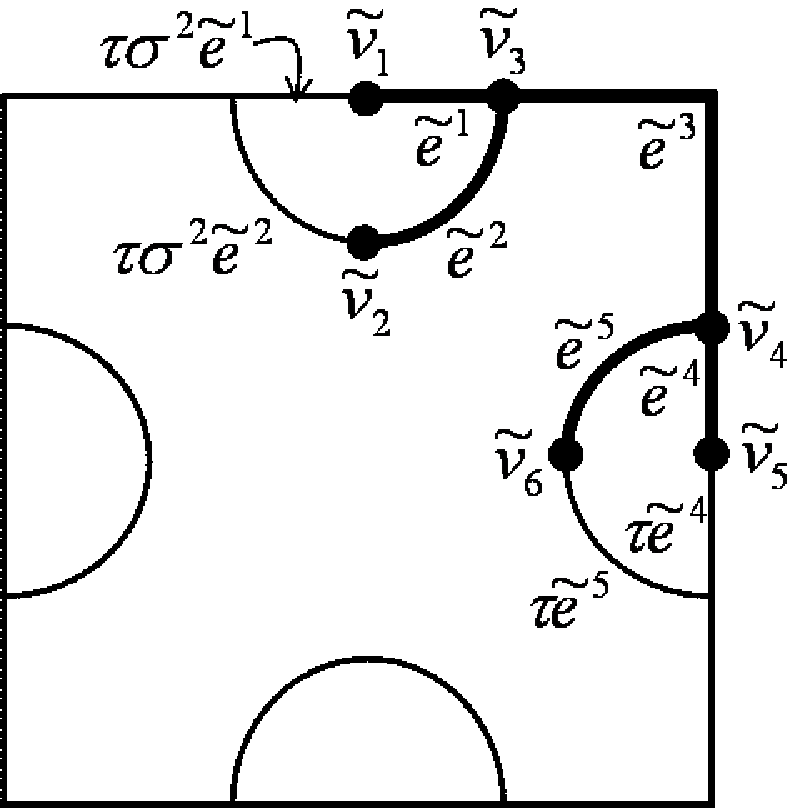} 
\par\end{center}

\begin{center}
(a) 
\par\end{center}%
\end{minipage}\hfill{}%
\begin{minipage}[c][1\totalheight][t]{0.35\columnwidth}%
\begin{center}
\includegraphics[width=0.6\textwidth]{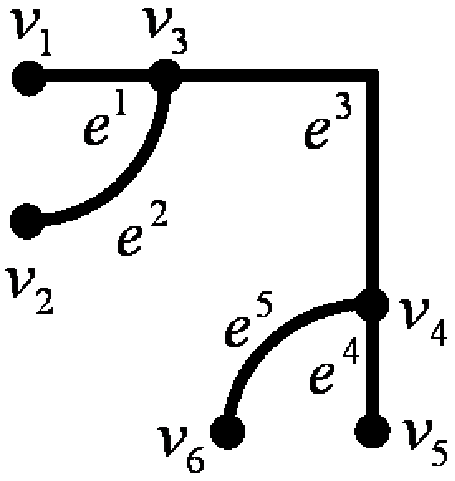} 
\par\end{center}

\medskip{}

\begin{center}
(b) 
\par\end{center}%
\end{minipage}\hfill{}

\caption{(a) The graph $\Gamma$ with the representatives of $\nicefrac{E}{H_{1}},\nicefrac{V}{H_{1}}$
marked in bold. (b) The resulting quotient $\nicefrac{\Gamma}{R_{1}}$.}

\label{fig:rigorous_quotient1} 
\end{figure}

\par\end{center}

First, we add {}``dummy'' vertices to the graph $\Gamma$ so that
no vertex is carried by the action of $H_{1}$ to one of its neighbors,
and choose representatives $\left\{ \tilde{e}^{i}\right\} _{i=1}^{5}$
for the orbits $\nicefrac{E}{H_{1}}$, and $\left\{ \tilde{v}_{k}\right\} _{k=1}^{6}$
for the orbits $\nicefrac{V}{H_{1}}$. These representatives are marked
in figure \ref{fig:rigorous_quotient1}(a) by bold lines and points.
The dummy vertices amongst the representatives are $\tilde{v}_{1},\tilde{v}_{2},\tilde{v}_{5},\tilde{v}_{6}$.
$R_{1}$ is one dimensional, and $d_{i}=1$ for all $i$ since the
stabilizers of all edges are trivial. Therefore, the quotient graph
is formed by taking one copy of each of the representative edges (figure
\ref{fig:rigorous_quotient1}(b)). Now, let us determine the boundary
conditions using (\ref{eq:final-boundary-conditions-1}), (\ref{eq:final-boundary-conditions-2}).
For all vertices we have $d=1$ and therefore $A_{\tilde{v}_{k}}\otimes I_{d}=A_{\tilde{v}_{k}}$
and $B_{\tilde{v}_{k}}\otimes I_{d}=B_{\tilde{v}_{k}}$. Consider
the vertex $v_{k}=v_{3}$ for which \begin{gather*}
n=3,\quad m=3,\quad d_{v_{3}}=3\\
g_{1}=e,\quad\nu_{1}=\mu_{1}=1,\qquad g_{2}=e,\quad\nu_{2}=\mu_{2}=2,\qquad g_{3}=e,\quad\nu_{3}=\mu_{3}=3\\
\mathfrak{G}=I_{3}\qquad\Theta=\left(\Theta'\otimes I_{d}\right)=\Theta'=I_{3}\end{gather*}
Plugging all this into (\ref{eq:final-boundary-conditions-1}), (\ref{eq:final-boundary-conditions-2})
and using the boundary conditions on $\tilde{v}_{3}$ which are given
by $A_{\tilde{v}_{3}}=\left(\begin{smallmatrix}1 & -1 & 0\\
0 & 1 & -1\\
0 & 0 & 0\end{smallmatrix}\right)$, $B_{\tilde{v}_{3}}=\left(\begin{smallmatrix}0 & 0 & 0\\
0 & 0 & 0\\
1 & 1 & 1\end{smallmatrix}\right)$ gives Neumann boundary conditions for $v_{3}$ as well: $A_{v_{3}}=A_{\tilde{v}_{3}},\quad B_{v_{3}}=B_{\tilde{v}_{3}}$.
Exactly the same treatment can be done for the vertex $v_{4}$ and
the same boundary conditions are obtained. The case is different for
the vertex $v_{5}$: \begin{gather*}
n=2,\quad m=1,\quad d_{v_{5}}=1\\
g_{1}=e,\quad\nu_{1}=\mu_{1}=4,\qquad g_{2}=\tau,\quad\nu_{2}=\mu_{1}=4\\
\mathfrak{G}=\left(\begin{smallmatrix}1 & 0\\
0 & -1\end{smallmatrix}\right)\qquad\Theta=\left(\Theta'\otimes I_{d}\right)=\Theta'=\left(\begin{smallmatrix}1\\
1\end{smallmatrix}\right)\end{gather*}
 The boundary conditions on $\tilde{v}_{5}$ are of Neumann type as
well: $A_{\tilde{v}_{5}}=\left(\begin{smallmatrix}1 & -1\\
0 & 0\end{smallmatrix}\right)$, $B_{\tilde{v}_{5}}=\left(\begin{smallmatrix}0 & 0\\
1 & 1\end{smallmatrix}\right)$. This time we obtain \begin{gather*}
A_{v_{5}}=\left(\begin{smallmatrix}1 & -1\\
0 & 0\end{smallmatrix}\right)\cdot\left(\begin{smallmatrix}1 & 0\\
0 & -1\end{smallmatrix}\right)\cdot\left(\begin{smallmatrix}1\\
1\end{smallmatrix}\right)=\left(\begin{smallmatrix}2\\
0\end{smallmatrix}\right),\\
B_{v_{5}}=\left(\begin{smallmatrix}0 & 0\\
1 & 1\end{smallmatrix}\right)\cdot\left(\begin{smallmatrix}1 & 0\\
0 & -1\end{smallmatrix}\right)\cdot\left(\begin{smallmatrix}1\\
1\end{smallmatrix}\right)=\left(\begin{smallmatrix}0\\
0\end{smallmatrix}\right).\end{gather*}
$A_{v_{5}}$ and $B_{v_{5}}$ are then reduced to square one dimensional
matrices as expected, by removing the second row in both of them.
We remain with $A_{v_{5}}=(2),\, B_{v_{5}}=(0)$ which means Dirichlet
boundary conditions on the vertex $v_{5}$. The same boundary conditions
are obtained for $v_{6}$. Similar derivation for vertices $v_{1},\, v_{2}$
gives Neumann boundary conditions for each one of them. The rigorous
construction thus gives us the same quotient graph that was obtained
by the intuitive method (figures \ref{fig:dihedral_triple}(a), \ref{fig:intuitive_quotient1}(b)).

\begin{figure}[!h]
\begin{centering}
\hfill{}%
\begin{minipage}[c][1\totalheight][t]{0.4\columnwidth}%
\begin{center}
\includegraphics[width=0.9\textwidth]{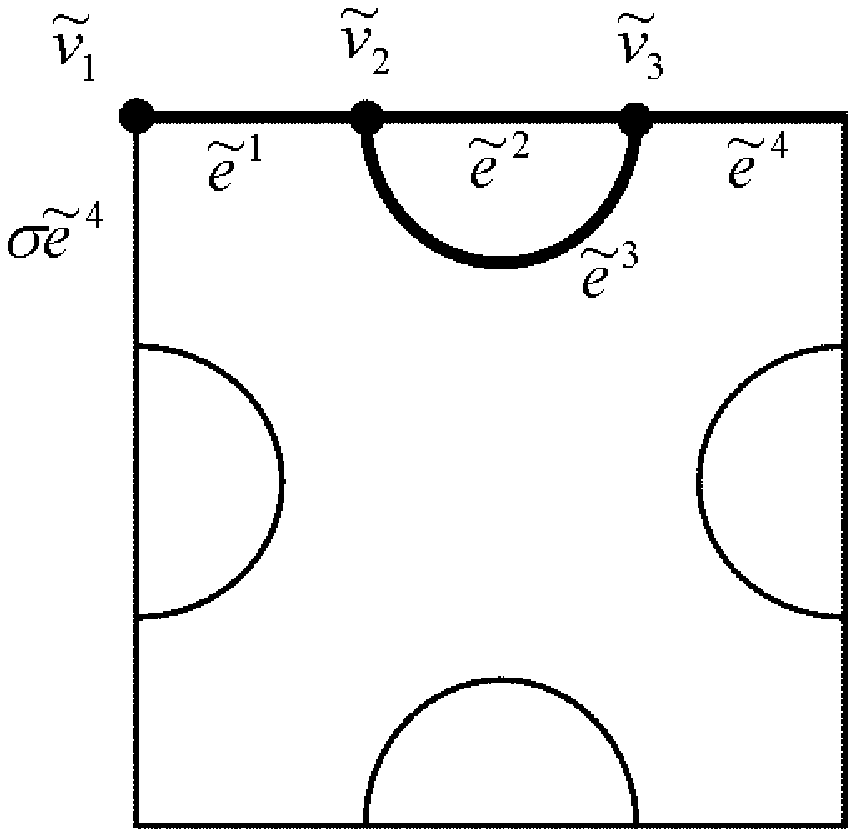} 
\par\end{center}

\begin{center}
(a) 
\par\end{center}%
\end{minipage}\hfill{}%
\begin{minipage}[t]{0.5\columnwidth}%
\begin{center}
\begin{minipage}[c][1\totalheight][t]{0.6\columnwidth}%
\begin{center}
\includegraphics[width=0.85\textwidth]{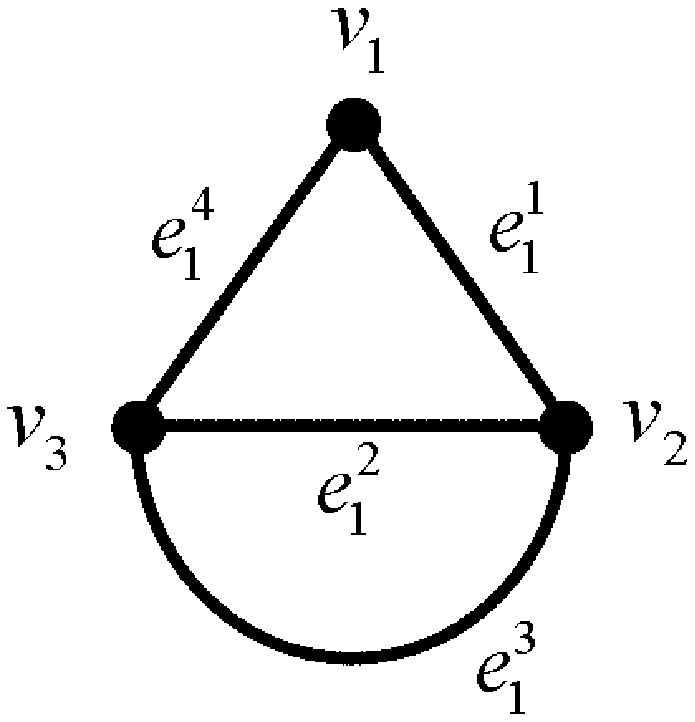} 
\par\end{center}

\medskip{}
\end{minipage}%
\begin{minipage}[c]{0.35\columnwidth}%
\begin{center}
$\begin{array}{l}
A_{v_{1}}=\left(\begin{smallmatrix}1 & i\\
0 & 0\end{smallmatrix}\right)\\
\\B_{v_{1}}=\left(\begin{smallmatrix}0 & 0\\
1 & -i\end{smallmatrix}\right)\end{array}$ 
\par\end{center}%
\end{minipage}
\par\end{center}

\begin{center}
(b) 
\par\end{center}%
\end{minipage}\hfill{}
\par\end{centering}

\caption{(a) The graph $\Gamma$ with the representatives of $\nicefrac{E}{H_{3}},\nicefrac{V}{H_{3}}$
marked in bold. (b) The resulting quotient $\nicefrac{\Gamma}{R_{3}}$.
$v_{2},\, v_{3}$ possess Neumann boundary conditions.}

\label{fig:quotient3} 
\end{figure}

The quotient $\nicefrac{\Gamma}{R_{2}}$ can be constructed in a similar
manner, and is shown in figure \ref{fig:dihedral_triple}(b). We proceed
to demonstrate the construction method for the quotient $\nicefrac{\Gamma}{R_{3}}$~%
\footnote{This result was obtained with G. Ben-Shach.%
}. We first add the corners of the square as dummy vertices to $\Gamma$
($\tilde{v}_{1}$ in figure \ref{fig:quotient3}(a) is one of them).
We are not obliged to do so, but it yields a quotient with simpler
boundary conditions. The choice of representatives for the edges and
the vertices is shown in figure \ref{fig:quotient3}(a) and the resulting
quotient in figure \ref{fig:quotient3}(b).\\
The vertices $v_{2}$ and $v_{3}$ have Neumann boundary conditions
exactly as their predecessors, $\tilde{v}_{2}$ and $\tilde{v}_{3}$.
For $v_{1}$ we obtain more interesting boundary conditions: \begin{gather*}
A_{\tilde{v}_{1}}=\left(\begin{smallmatrix}1 & -1\\
0 & 0\end{smallmatrix}\right),\quad B_{\tilde{v}_{1}}=\left(\begin{smallmatrix}0 & 0\\
1 & 1\end{smallmatrix}\right)\\
n=2,\quad m=2,\quad d_{v_{1}}=2,\\
g_{1}=e,\quad\nu_{1}=\mu_{1}=1,\qquad g_{2}=\sigma,\quad\nu_{2}=\mu_{2}=4\\
\mathfrak{G}=\left(\begin{smallmatrix}1 & 0\\
0 & -i\end{smallmatrix}\right),\qquad\Theta=\left(\Theta'\otimes I_{d}\right)=\Theta'=I_{2}\:,\end{gather*}
which gives \begin{equation}
A_{v_{1}}=\left(\begin{smallmatrix}1 & i\\
0 & 0\end{smallmatrix}\right),\qquad B_{v_{1}}=\left(\begin{smallmatrix}0 & 0\\
1 & -i\end{smallmatrix}\right)\:.\label{eq:tripleH3_v1_cond}\end{equation}
Non-formally speaking, the vertex $v_{1}$ {}``applies a factor of
$i$'' to the functions that cross it. The resulting graph is the
one that was shown in figure \ref{fig:dihedral_triple}(c).

\medskip{}

In order to exhaust this example, we observe that $\ind_{H_{1}}^{G}R_{1}\cong\ind_{H_{2}}^{G}R_{2}\cong\ind_{H_{3}}^{G}R_{3}$
is the two-dimensional irreducible representation of $D_{4}$, which
we denote by $R$. By theorem \ref{thm:mainthm}, the isospectral
family of the three graphs given in figure \ref{fig:dihedral_triple}
can be extended by adding any graph which is $\nicefrac{\Gamma}{R}$.
We therefore construct now such a graph. Let us use the intuitive
approach first. Recall that \eqref{eq:hom-simple-isotip} was the
key for the intuitive construction of $\nicefrac{\Gamma}{R_{1}}$.
Analogously to \eqref{eq:hom-simple-isotip}, we make the observation
that encoding $\mathrm{Hom}_{\mathbb{C}G}\left(R,\mathcal{H}\left(\Gamma\right)\right)$,
the r.h.s. of \eqref{eq:quotient-definition}, is similar in nature
to encoding $\left(\mathcal{H}\left(\Gamma\right)\right)^{R}$, the
$R$-isotypic component of $\mathcal{H}\left(\Gamma\right)$, as due
to the simplicity of $R$ as a $\mathbb{C}G$-module the two are isomorphic.
This can be understood as follows: making a choice of a basis $\left\{ b_{1},b_{2}\right\} $
for $R$, and given a function $\tilde{f}\in\mathrm{Hom}_{\mathbb{C}G}\left(R,\mathcal{H}\left(\Gamma\right)\right)$,
we have that $\tilde{f}_{b_{1}},\,\tilde{f}_{b_{2}}\in\left(\mathcal{H}\left(\Gamma\right)\right)^{R}$
and furthermore $\left\{ \tilde{f}_{b_{1}},\,\tilde{f}_{b_{2}}\right\} $
spans over $\mathbb{C}$ a $\mathbb{C}G$-module isomorphic to $R$.
In order to exhibit the general behavior we avoid sparse matrices,
and pick a basis $\left\{ b_{1},b_{2}\right\} $ for which the matrix
representation of $R$ is \begin{equation}
\left\{ \begin{array}{ll}
\tau\sigma^{2}\mapsto\frac{1}{2}\left(\begin{array}{cc}
-1 & -\sqrt{3}\\
-\sqrt{3} & 1\end{array}\right), & \tau\sigma^{3}\mapsto\frac{1}{2}\left(\begin{array}{cc}
\sqrt{3} & -1\\
-1 & -\sqrt{3}\end{array}\right)\end{array}\right\} \,.\label{eq:matrix_rep}\end{equation}
It is enough to consider only the matrices of these two elements for
the construction of the quotient.

\begin{center}
\begin{figure}[!h]
\begin{minipage}[c]{0.64\columnwidth}%
\begin{center}
\includegraphics[width=0.95\textwidth]{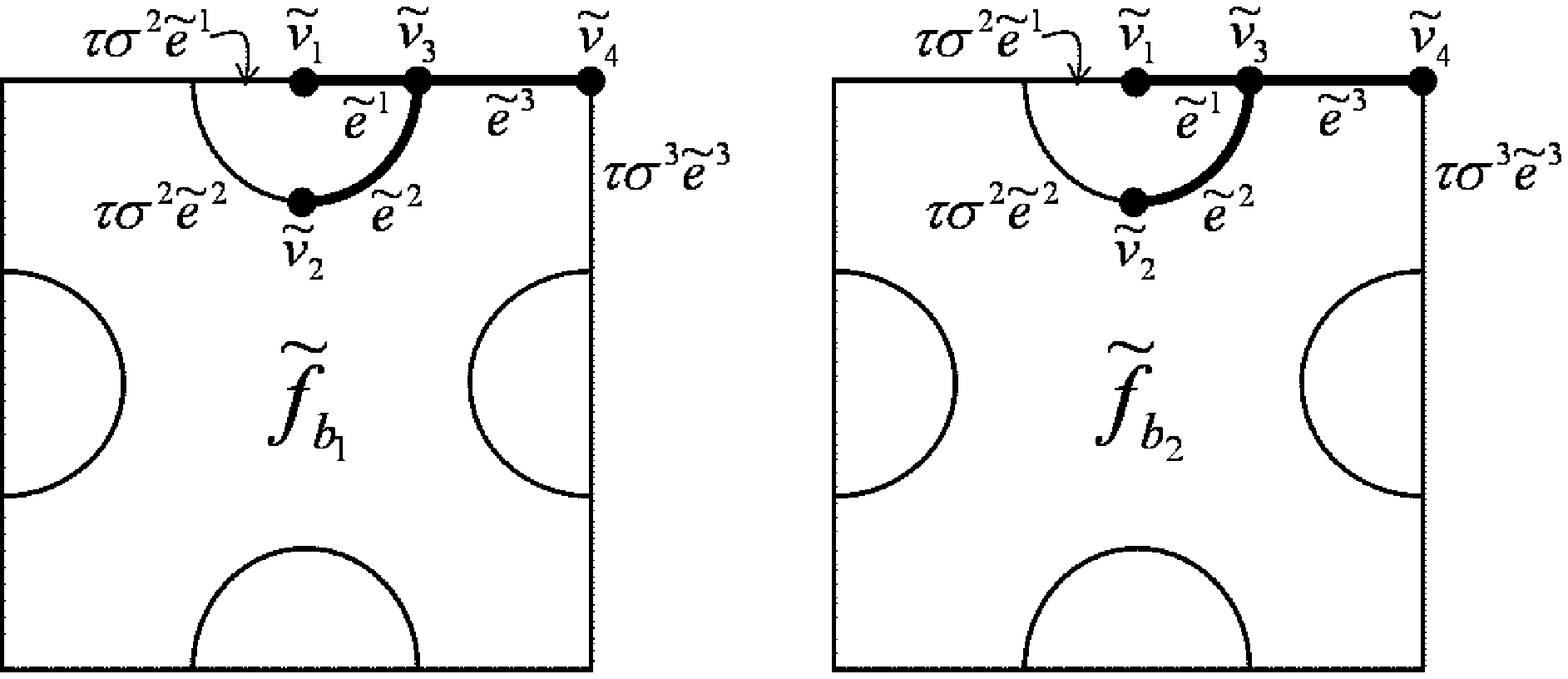} 
\par\end{center}

\begin{center}
(a) 
\par\end{center}%
\end{minipage}\hfill{}%
\begin{minipage}[c]{0.31\columnwidth}%
\begin{center}
\includegraphics[width=0.85\textwidth]{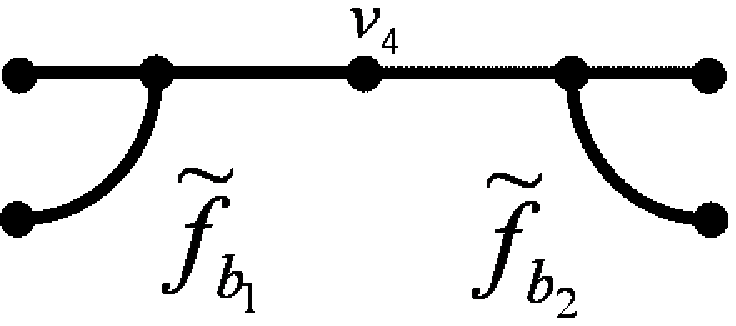} 
\par\end{center}

\medskip{}

\begin{center}
(b) 
\par\end{center}%
\end{minipage}

\caption{(a) Two copies of the graph $\Gamma$ with the representatives of
$\nicefrac{E}{D_{4}},\nicefrac{V}{D_{4}}$ marked in bold. These two
copies are merely a visualization of the {}``basis functions'' $\tilde{f}_{b_{1}},\tilde{f}_{b_{2}}$
on $\Gamma$. (b) The first stage in the formation of $\nicefrac{\Gamma}{R}$
is the gluing of both copies in the vertex $v_{4}$, with the boundary
conditions given in \eqref{eq:v4_values}, \eqref{eq:v4_derivatives}
.}

\label{fig:quotient4ab} 
\end{figure}

\par\end{center}

\begin{figure}[!h]
\begin{centering}
\includegraphics[scale=0.45]{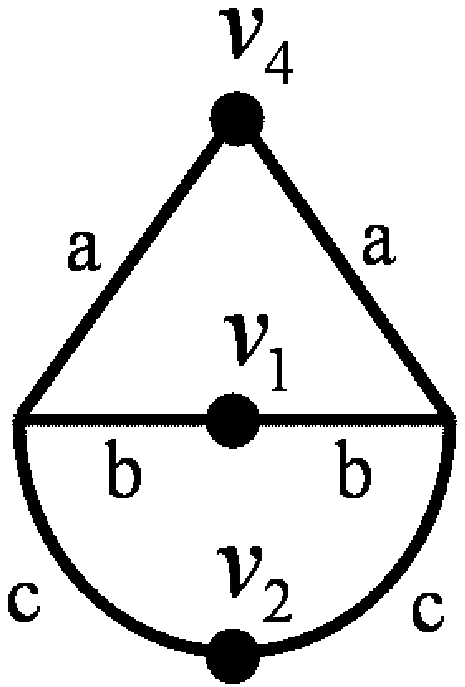} 
\par\end{centering}

\caption{The quotient graph $\nicefrac{\Gamma}{R}$ which is isospectral to
the graphs in figure \ref{fig:dihedral_triple}. The boundary conditions
are as described in \eqref{eq:v4_bc}, \eqref{eq:v12_bc}.}

\label{fig:quotient4c} 
\end{figure}

Examine the properties of $\tilde{f}_{b_{1}},\,\tilde{f}_{b_{2}}$
that follow from the above matrix representation (figure \ref{fig:quotient4ab}(a)).
Since $\tilde{f}\in\mathrm{Hom}_{\mathbb{C}G}\left(R,\mathcal{H}\left(\Gamma\right)\right)$
we have $\tau\sigma^{3}\tilde{f}_{b_{1}}=\tilde{f}_{\left(\tau\sigma^{3}\right)^{-1}b_{1}}=\tilde{f}_{\tau\sigma^{3}b_{1}}$,
and thus the first column of the matrix representing $\tau\sigma^{3}$
tells us that \begin{eqnarray}
\tau\sigma^{3}\tilde{f}_{b_{1}} & = & \nicefrac{\sqrt{3}}{2}\tilde{f}_{b_{1}}-\nicefrac{1}{2}\tilde{f}_{b_{2}}\label{eq:diagonal_reflection_values}\\
\tau\sigma^{3}\tilde{f}'_{b_{1}} & = & \nicefrac{\sqrt{3}}{2}\tilde{f}'_{b_{1}}-\nicefrac{1}{2}\tilde{f}'_{b_{2}}\label{eq:diagonal_reflection_deriv}\end{eqnarray}
and enables us to relate the values and the derivatives of $\tilde{f}_{b_{1}},\tilde{f}_{b_{2}}$
on the vertex $\tilde{v}_{4}$. Since $\tilde{v}_{4}$ is a fixed
point under the action of $\tau\sigma^{3}$ and there are Neumann
boundary conditions on it, we have that \begin{eqnarray}
\left(\tau\sigma^{3}\tilde{f}_{b_{1}}\right)\at_{\tilde{e}_{3}}(\tilde{v}_{4}) & = & \tilde{f}_{b_{1}}\at_{\tilde{e}_{3}}(\tilde{v}_{4})\label{eq:v4_temp_values}\\
\left(\tau\sigma^{3}\tilde{f}'_{b_{1}}\right)\at_{\tilde{e}_{3}}(\tilde{v}_{4}) & = & -\tilde{f}'_{b_{1}}\at_{\tilde{e}_{3}}(\tilde{v}_{4})\,.\label{eq:v4_temp_deriv}\end{eqnarray}
Evaluating \eqref{eq:diagonal_reflection_values} on $v_{4}$ and
combining this with \eqref{eq:v4_temp_values} gives \begin{equation}
\left(1-\nicefrac{\sqrt{3}}{2}\right)\tilde{f}_{b_{1}}\at_{\tilde{e}_{3}}(\tilde{v}_{4})+\nicefrac{1}{2}\tilde{f}_{b_{2}}\at_{\tilde{e}_{3}}(\tilde{v}_{4})=0\,.\label{eq:v4_values}\end{equation}
Similarly, from \eqref{eq:diagonal_reflection_deriv} and \eqref{eq:v4_temp_deriv}
we obtain \begin{equation}
\left(-1-\nicefrac{\sqrt{3}}{2}\right)\tilde{f}_{b_{1}}'\at_{\tilde{e}_{3}}(\tilde{v}_{4})+\nicefrac{1}{2}\tilde{f}_{b_{2}}'\at_{\tilde{e}_{3}}(\tilde{v}_{4})=0\,.\label{eq:v4_derivatives}\end{equation}
We may therefore think of two copies of the graphs. Each of the basis
functions $\tilde{f}_{b_{1}},\,\tilde{f}_{b_{2}}$ resides on one
of the copies, and the relations between the values and the derivatives
of the functions allow us to take a subgraph out of each copy (marked
in bold in figure \ref{fig:quotient4ab}(a)) and glue both of them
together with the appropriate boundary conditions. The first stage
in this gluing process, visualized in figure \ref{fig:quotient4ab}(b),
is to identify the vertex $\tilde{v}_{4}$ in the two copies and turn
it into the vertex $v_{4}$ of the quotient with the boundary conditions
that were derived in \eqref{eq:v4_values}, \eqref{eq:v4_derivatives}:
\begin{eqnarray}
A_{v_{4}}=\left(\begin{smallmatrix}1-\nicefrac{\sqrt{3}}{2} & \nicefrac{1}{2}\\
0 & 0\end{smallmatrix}\right) & \quad,\quad & B_{v_{4}}=\left(\begin{smallmatrix}0 & 0\\
-1-\nicefrac{\sqrt{3}}{2} & \nicefrac{1}{2}\end{smallmatrix}\right)\,.\label{eq:v4_bc}\end{eqnarray}
After treating similarly vertices $\tilde{v}_{1},\tilde{v}_{2}$ we
get the quotient $\nicefrac{\Gamma}{R}$ (figure \ref{fig:quotient4c})
whose remaining boundary conditions are given by: \begin{gather}
A_{v_{1}}=A_{v_{2}}=\left(\begin{smallmatrix}\nicefrac{3}{2} & \nicefrac{\sqrt{3}}{2}\\
0 & 0\end{smallmatrix}\right)\quad,\quad B_{v_{1}}=B_{v_{2}}=\left(\begin{smallmatrix}0 & 0\\
-\nicefrac{1}{2} & \nicefrac{\sqrt{3}}{2}\end{smallmatrix}\right)\,.\label{eq:v12_bc}\end{gather}

\medskip{}

\begin{figure}[!h]
\hfill{}%
\begin{minipage}[t]{0.35\columnwidth}%
\begin{center}
\includegraphics[width=0.8\textwidth]{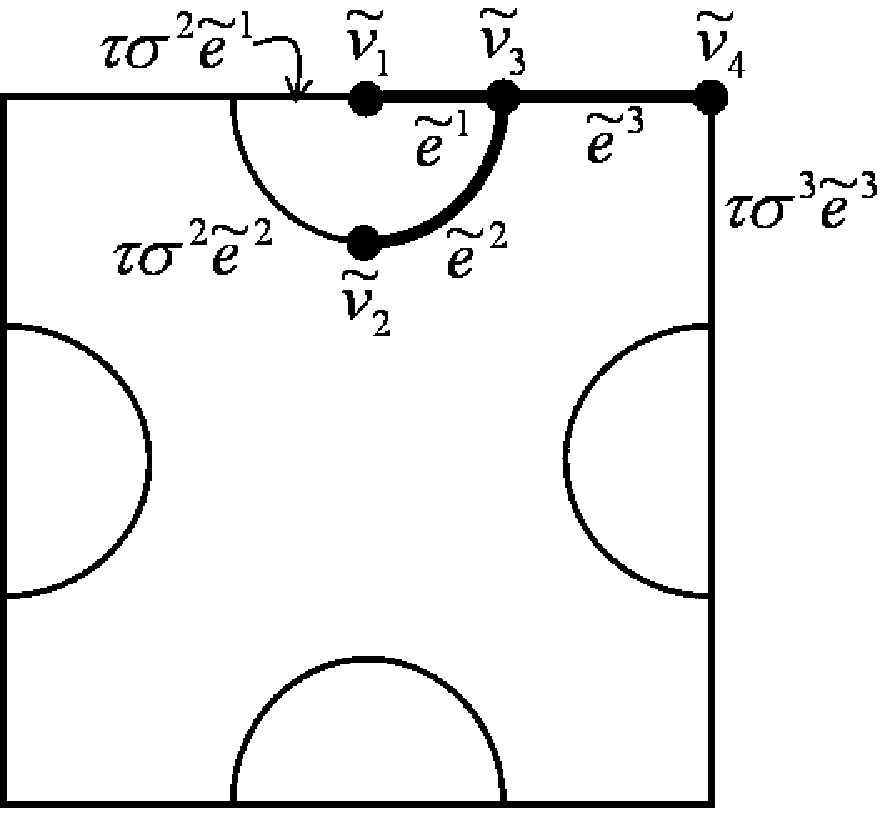} 
\par\end{center}

\begin{center}
(a) 
\par\end{center}%
\end{minipage}\hfill{}%
\begin{minipage}[t]{0.45\columnwidth}%
\begin{center}
\includegraphics[width=0.95\textwidth]{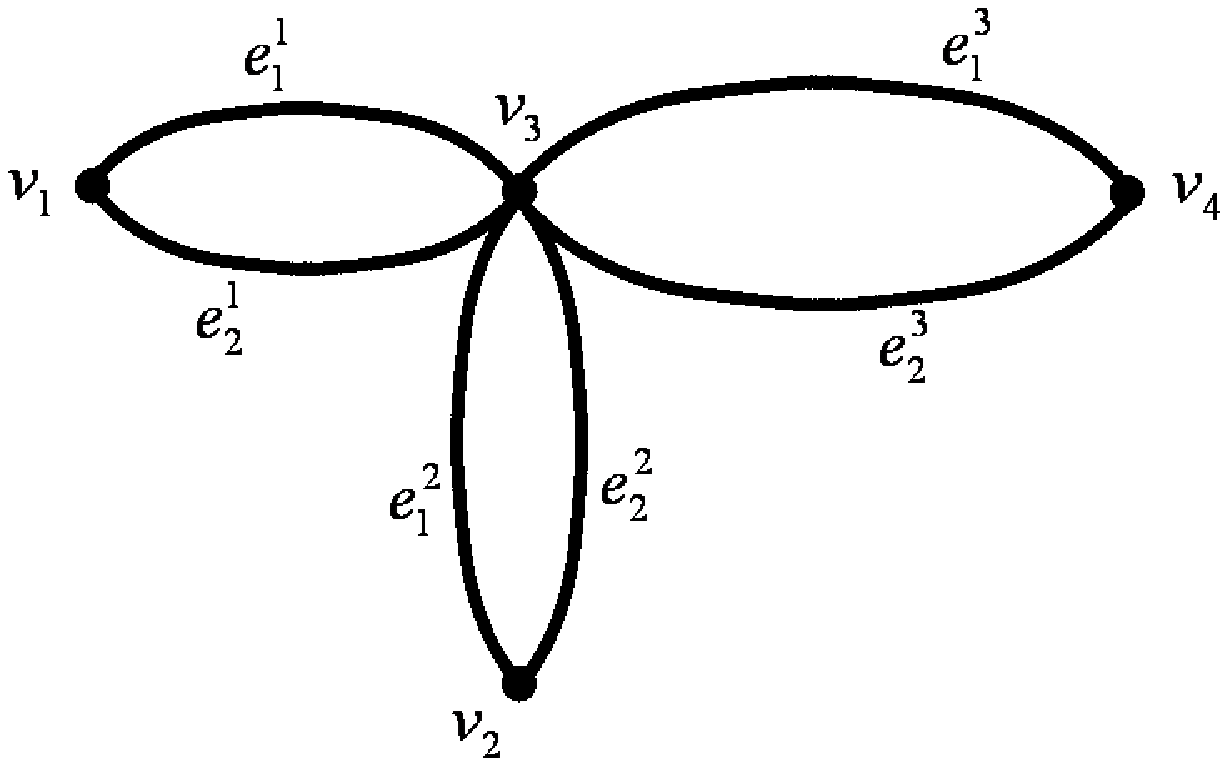} 
\par\end{center}

\begin{center}
(b) 
\par\end{center}%
\end{minipage}\hfill{}

\caption{(a) The graph $\Gamma$ with the representatives of $\nicefrac{E}{D_{4}},\nicefrac{V}{D_{4}}$
marked in bold. (b) The resulting quotient $\nicefrac{\Gamma}{R}$.}

\label{fig:quotient4de} 
\end{figure}

We now use the rigorous approach for the same quotient, $\nicefrac{\Gamma}{R}$.
The representatives of the orbits $\nicefrac{E}{G}$ are $\left\{ \tilde{e}^{i}\right\} _{i=1}^{3}$
and the representatives of $\nicefrac{V}{G}$ are $\left\{ \tilde{v}_{k}\right\} _{k=1}^{4}$
(figure \ref{fig:quotient4de}(a)). This time the representation is
not one dimensional ($d=2$) so there are additional details to consider.
First, note that we have two copies of each representative of $\nicefrac{E}{G}$
in the quotient and both of the copies survive since all edges have
trivial stabilizers (figure \ref{fig:quotient4de}(b)). This last
observation ensures that we can take $B^{i}=B$ for all $i$ (i.e.,
the same basis for all edges). We again take $B$ to be the basis
for which the matrix representation of $R$ is \eqref{eq:matrix_rep}.

We treat the boundary conditions at the vertices one by one:
\begin{itemize}
\item $v_{4}$ has the following data: \begin{gather*}
n=2,\quad m=1,\quad d_{v_{4}}=2,\qquad g_{1}=e,\quad\nu_{1}=\mu_{1}=3,\qquad g_{2}=\tau\sigma^{3},\quad\nu_{2}=\mu_{1}=3\\
\Theta=\left(\Theta'\otimes I_{d}\right)=\left(\begin{smallmatrix}1\\
1\end{smallmatrix}\right)\otimes I_{2}=\left(\begin{smallmatrix}1 & 0\\
0 & 1\\
1 & 0\\
0 & 1\end{smallmatrix}\right)\end{gather*}
 $A_{\tilde{v_{4}}},B_{\tilde{v_{4}}}$ are the regular Neumann matrices
and we therefore obtain \begin{eqnarray*}
A_{v_{4}}= & \left(\begin{smallmatrix}1 & 0 & -1 & 0\\
0 & 1 & 0 & -1\\
0 & 0 & 0 & 0\\
0 & 0 & 0 & 0\end{smallmatrix}\right)\cdot\left(\begin{smallmatrix}1 & 0 & 0 & 0\\
0 & 1 & 0 & 0\\
0 & 0 & \nicefrac{\sqrt{3}}{2} & -\nicefrac{1}{2}\\
0 & 0 & -\nicefrac{1}{2} & -\nicefrac{\sqrt{3}}{2}\end{smallmatrix}\right)\ \cdot\left(\begin{smallmatrix}1 & 0\\
0 & 1\\
1 & 0\\
0 & 1\end{smallmatrix}\right) & =\left(\begin{smallmatrix}1-\nicefrac{\sqrt{3}}{2} & \nicefrac{1}{2}\\
\nicefrac{1}{2} & 1+\nicefrac{\sqrt{3}}{2}\\
0 & 0\\
0 & 0\end{smallmatrix}\right)\\
B_{v_{4}}= & \left(\begin{smallmatrix}0 & 0 & 0 & 0\\
0 & 0 & 0 & 0\\
1 & 0 & 1 & 0\\
0 & 1 & 0 & 1\end{smallmatrix}\right)\cdot\left(\begin{smallmatrix}1 & 0 & 0 & 0\\
0 & 1 & 0 & 0\\
0 & 0 & \nicefrac{\sqrt{3}}{2} & -\nicefrac{1}{2}\\
0 & 0 & -\nicefrac{1}{2} & -\nicefrac{\sqrt{3}}{2}\end{smallmatrix}\right)\ \cdot\left(\begin{smallmatrix}1 & 0\\
0 & 1\\
1 & 0\\
0 & 1\end{smallmatrix}\right) & =\left(\begin{smallmatrix}0 & 0\\
0 & 0\\
1+\nicefrac{\sqrt{3}}{2} & -\nicefrac{1}{2}\\
-\nicefrac{1}{2} & 1-\nicefrac{\sqrt{3}}{2}\end{smallmatrix}\right)\end{eqnarray*}
 Noting that both $A_{v_{4}}$ and $B_{v_{4}}$ are of rank one, we
see that they express the same boundary conditions as given in \eqref{eq:v4_bc}. 
\item $v_{1}$ obviously has the same boundary conditions as $v_{2}$. We
examine $v_{1}$: \begin{gather*}
n=2,\quad m=1,\quad d_{v_{1}}=2,\qquad g_{1}=e,\quad\nu_{1}=\mu_{1}=1,\qquad g_{2}=\tau\sigma^{2},\quad\nu_{2}=\mu_{1}=1\\
\Theta=\left(\Theta'\otimes I_{d}\right)=\left(\begin{smallmatrix}1\\
1\end{smallmatrix}\right)\otimes I_{2}=\left(\begin{smallmatrix}1 & 0\\
0 & 1\\
1 & 0\\
0 & 1\end{smallmatrix}\right)\end{gather*}
 Again, $A_{\tilde{v}_{1}}$ and $B_{\tilde{v}_{1}}$ are the regular
Neumann matrices and we get: \begin{eqnarray*}
A_{v_{1}}= & \left(\begin{smallmatrix}1 & 0 & -1 & 0\\
0 & 1 & 0 & -1\\
0 & 0 & 0 & 0\\
0 & 0 & 0 & 0\end{smallmatrix}\right)\cdot\left(\begin{smallmatrix}1 & 0 & 0 & 0\\
0 & 1 & 0 & 0\\
0 & 0 & -\nicefrac{1}{2} & -\nicefrac{\sqrt{3}}{2}\\
0 & 0 & -\nicefrac{\sqrt{3}}{2} & \nicefrac{1}{2}\end{smallmatrix}\right)\ \cdot\left(\begin{smallmatrix}1 & 0\\
0 & 1\\
1 & 0\\
0 & 1\end{smallmatrix}\right) & =\left(\begin{smallmatrix}\nicefrac{3}{2} & \nicefrac{\sqrt{3}}{2}\\
\nicefrac{\sqrt{3}}{2} & \nicefrac{1}{2}\\
0 & 0\\
0 & 0\end{smallmatrix}\right)\\
B_{v_{1}}= & \left(\begin{smallmatrix}0 & 0 & 0 & 0\\
0 & 0 & 0 & 0\\
1 & 0 & 1 & 0\\
0 & 1 & 0 & 1\end{smallmatrix}\right)\cdot\left(\begin{smallmatrix}1 & 0 & 0 & 0\\
0 & 1 & 0 & 0\\
0 & 0 & -\nicefrac{1}{2} & -\nicefrac{\sqrt{3}}{2}\\
0 & 0 & -\nicefrac{\sqrt{3}}{2} & \nicefrac{1}{2}\end{smallmatrix}\right)\ \cdot\left(\begin{smallmatrix}1 & 0\\
0 & 1\\
1 & 0\\
0 & 1\end{smallmatrix}\right) & =\left(\begin{smallmatrix}0 & 0\\
0 & 0\\
\nicefrac{1}{2} & -\nicefrac{\sqrt{3}}{2}\\
-\nicefrac{\sqrt{3}}{2} & \nicefrac{3}{2}\end{smallmatrix}\right)\end{eqnarray*}
 which are matrices of rank one and again we may reduce these matrices
into two dimensional ones which are exactly those given in \eqref{eq:v12_bc}. 
\item The case of $v_{3}$ is a bit more interesting: \begin{gather*}
n=3,\quad m=3,\quad d_{v_{3}}=6,\\
g_{1}=e,\quad\nu_{1}=\mu_{1}=1,\qquad g_{2}=e,\quad\nu_{2}=\mu_{2}=2,\qquad g_{3}=e,\quad\nu_{3}=\mu_{3}=3\\
\Theta=\left(\Theta'\otimes I_{d}\right)=I_{3}\otimes I_{2}=I_{6}\end{gather*}
As $A_{\tilde{v}_{3}}$ and $B_{\tilde{v}_{3}}$ are Neumann matrices,
we have \begin{eqnarray*}
A_{v_{3}} & = & \left(\begin{smallmatrix}1 & 0 & -1 & 0 & 0 & 0\\
0 & 1 & 0 & -1 & 0 & 0\\
0 & 0 & 1 & 0 & -1 & 0\\
0 & 0 & 0 & 1 & 0 & -1\\
0 & 0 & 0 & 0 & 0 & 0\\
0 & 0 & 0 & 0 & 0 & 0\end{smallmatrix}\right)\cdot I_{6}\cdot I_{6}=\left(\begin{smallmatrix}1 & 0 & -1 & 0 & 0 & 0\\
0 & 1 & 0 & -1 & 0 & 0\\
0 & 0 & 1 & 0 & -1 & 0\\
0 & 0 & 0 & 1 & 0 & -1\\
0 & 0 & 0 & 0 & 0 & 0\\
0 & 0 & 0 & 0 & 0 & 0\end{smallmatrix}\right)\\
B_{v_{3}} & = & \left(\begin{smallmatrix}0 & 0 & 0 & 0 & 0 & 0\\
0 & 0 & 0 & 0 & 0 & 0\\
0 & 0 & 0 & 0 & 0 & 0\\
0 & 0 & 0 & 0 & 0 & 0\\
1 & 0 & 1 & 0 & 1 & 0\\
0 & 1 & 0 & 1 & 0 & 1\end{smallmatrix}\right)\cdot I_{6}\cdot I_{6}=\left(\begin{smallmatrix}0 & 0 & 0 & 0 & 0 & 0\\
0 & 0 & 0 & 0 & 0 & 0\\
0 & 0 & 0 & 0 & 0 & 0\\
0 & 0 & 0 & 0 & 0 & 0\\
1 & 0 & 1 & 0 & 1 & 0\\
0 & 1 & 0 & 1 & 0 & 1\end{smallmatrix}\right)\:,\end{eqnarray*}
and we see that the above boundary conditions separate the edges into
two sets, $\left\{ e_{1}^{1},e_{1}^{2},e_{1}^{3}\right\} $ and $\left\{ e_{2}^{1},e_{2}^{2},e_{2}^{3}\right\} $,
each dominated by a regular Neumann condition. This enables us to
split the vertex $v_{3}$ into two distinct vertices of degree $3$,
each connected to a different set of edges and possessing Neumann
boundary conditions. We remark that this would happen for any choice
of basis for $R$, as here $g_{1}=g_{2}=g_{3}=e$.
\end{itemize}
Note that the resulting quotient is the same as was obtained previously
(figure \ref{fig:quotient4c}).

\medskip{}

Finally, we repeat the construction for an arbitrary choice of basis
which yields an orthogonal matrix representation for $R$. We can
parametrize such a representation in the following way: \[
\left\{ \begin{array}{l}
\tau\sigma^{2}\mapsto\left(\begin{array}{cc}
\cos^{2}\theta-\sin^{2}\theta & -2\cos\theta\sin\theta\\
-2\cos{\theta}\sin{\theta} & -\cos^{2}{\theta}+\sin^{2}{\theta}\end{array}\right),\\
\tau\sigma^{3}\mapsto\left(\begin{array}{cc}
2\cos{\theta}\sin{\theta} & \cos^{2}{\theta}-\sin^{2}{\theta}\\
\cos^{2}{\theta}-\sin^{2}{\theta} & -2\cos{\theta}\sin{\theta}\end{array}\right)\end{array}\right\} .\]
 For example, the basis we chose in \eqref{eq:matrix_rep} is obtained
by $\theta=\nicefrac{\pi}{3}$. As remarked, $v_{3}$ always splits
into two vertices with Neumann conditions, so that figure \ref{fig:quotient4c}
can describe the quotient with respect to any basis. For the parametrization
above, we obtain the following boundary conditions: \begin{align*}
A_{v_{1}}=A_{v_{2}}=\left(\begin{smallmatrix}2\sin^{2}{\theta} & \sin2{\theta}\\
\sin2{\theta} & 2-2\sin^{2}{\theta}\\
0 & 0\\
0 & 0\end{smallmatrix}\right) & \qquad & A_{v_{4}}=\left(\begin{smallmatrix}1-\sin2{\theta} & 2\sin^{2}{\theta}-1\\
2\sin^{2}{\theta}-1 & 1+\sin2{\theta}\\
0 & 0\\
0 & 0\end{smallmatrix}\right)\\
B_{v_{1}}=B_{v_{2}}=\left(\begin{smallmatrix}0 & 0\\
0 & 0\\
2-2\sin^{2}{\theta} & -\sin2{\theta}\\
-\sin2{\theta} & 2\sin^{2}{\theta}\end{smallmatrix}\right) & \qquad & B_{v_{4}}=\left(\begin{smallmatrix}0 & 0\\
0 & 0\\
1+\sin2{\theta} & 1-2\sin^{2}{\theta}\\
1-2\sin^{2}{\theta} & 1-\sin2{\theta}\end{smallmatrix}\right)\end{align*}
All of these matrices are of rank one, and can therefore be reduced
to square ones by deleting the appropriate rows%
\footnote{However, there is no a priori reduction which is valid for all $\theta$!%
}. We thus get a continuous family of isospectral graphs. Examine two
members of this family: $\theta=0$ and $\theta=\nicefrac{3\pi}{4}$.
The boundary conditions for the case $\theta=0$ are: \begin{align*}
A_{v_{1}}=A_{v_{2}}=\left(\begin{smallmatrix}0 & 2\\
0 & 0\end{smallmatrix}\right)\qquad & \qquad A_{v_{4}}=\left(\begin{smallmatrix}1 & -1\\
0 & 0\end{smallmatrix}\right)\\
B_{v_{1}}=B_{v_{2}}=\left(\begin{smallmatrix}2 & 0\\
0 & 0\end{smallmatrix}\right)\qquad & \qquad B_{v_{4}}=\left(\begin{smallmatrix}0 & 0\\
1 & 1\end{smallmatrix}\right)\end{align*}
 When applying this to figure \ref{fig:quotient4c}, we notice that
the vertices $v_{1}$, $v_{2}$ do not stay vertices of degree two,
but rather, each of them splits into two vertices of degree one, one
with Dirichlet boundary condition, and the other with Neumann. The
vertex $v_{4}$, however, stays connected and obtains Neumann boundary
conditions. Observe that the resulting quotient is the one that we
have already obtained as $\nicefrac{\Gamma}{R_{1}}$ (figure \ref{fig:dihedral_triple}(a)).
In a similar manner, the quotient $\nicefrac{\Gamma}{R_{2}}$ (figure
\ref{fig:dihedral_triple}(b)) is obtained from the choice $\theta=\nicefrac{3\pi}{4}$.
We conclude by pointing out that the graph described in figure \ref{fig:quotient4c}
is a good prototype for the mentioned isospectral family, yet it might
also be misleading, since there are members of the family whose boundary
conditions tear apart the edges connected to some of the vertices
and thus change the connectivity of the graph. One should also pay
attention to the fact that we have treated only orthogonal representations
of $D_{4}$. These are not the most general ones, and we may extend
the isospectral family presented above by considering the broader
case of all matrix representations of $R$. In particular, the quotient
$\nicefrac{\Gamma}{R_{3}}$ (figure \ref{fig:dihedral_triple}(c))
is obtained from the unitary representation\[
\left\{ \begin{array}{cc}
\sigma\mapsto\left(\begin{array}{cc}
i & 0\\
0 & -i\end{array}\right), & \tau\mapsto\left(\begin{array}{cc}
0 & -1\\
-1 & 0\end{array}\right)\end{array}\right\} \:.\]

\section{Isospectral manifolds and stratifolds}

\label{sec:manifolds}

If $\Gamma$ is a Riemannian manifold equipped with an action of a
finite group $G$, then $C^{^{\infty}}\negthickspace\left(\Gamma\right)$
is again a module over $\mathbb{C}G\left[x\right]$, with $x$ acting
as the Laplace-Beltrami operator $\Delta$. If however $\Gamma$ has
a boundary, at which differential boundary conditions are imposed
on $C^{^{\infty}}\negthickspace\left(\Gamma\right)$, then in general
it is no longer closed under $\Delta$. In order to treat this case
as well, we limit our attention to the subspace of $C^{^{\infty}}\negthickspace\left(\Gamma\right)$
spanned by $\Delta$'s eigenfunctions, which we again denote by $\mathcal{H}\left(\Gamma\right)$.
Assuming that the boundary conditions are linear, $\mathcal{H}\left(\Gamma\right)$
is closed under $\Delta$ and is therefore a $\mathbb{C}G\left[x\right]$-module
as before. Section \ref{sec:algebra} is naturally generalized to
these settings:
\begin{itemize}
\item For a representation $R$ of $G$, we define a $\nicefrac{\Gamma}{R}$-manifold
to be a Riemannian manifold (possibly with boundary, at which homogeneous
conditions are imposed) $\Gamma'$, such that there is an isomorphism\begin{equation}
\mathcal{H}\left(\Gamma'\right)\cong\mathrm{Hom}_{\mathbb{C}G}\left(R,\mathcal{H}\left(\Gamma\right)\right)\label{eq:manifold-quotient-def}\end{equation}
 intertwining the Laplace-Beltrami operator.
\item For a representation $R$ of $H\leq G$, $\nicefrac{\Gamma}{R}$ and
$\nicefrac{\Gamma}{\ind_{H}^{G}R}$ are isospectral; therefore, $\nicefrac{\Gamma}{\mathbb{C}G}$
is isospectral to $\Gamma$, and for representations $R_{1},R_{2}$
of $H_{1},H_{2}\leq G$ satisfying $\ind_{H_{1}}^{G}R_{1}\cong\ind_{H_{2}}^{G}R_{2}$,
$\nicefrac{\Gamma}{R_{1}}$ and $\nicefrac{\Gamma}{R_{2}}$ are isospectral. 
\end{itemize}
The main advantage of quantum graphs for our purposes is that under
fairly moderate assumptions (e.g., self-adjoint Laplacian or a free
action) one can build a quotient for every representation, as is demonstrated
in section \ref{sec:quotient graph}. 

Graphs are one-dimensional manifolds with singularities (at the vertices),
and it is these singularities that we exploit, by endowing them with
the appropriate boundary conditions, to encapsulate the restrictions
arising from a choice of a representation. In higher dimensions, manifolds
with a boundary, carrying Neumann, Dirichlet, or a more complicated
boundary condition, are a generalization of this idea, and one goal
of this section is to demonstrate that some known isospectral examples
of such objects can be understood using our theory. That is, we show
that for some known isospectral pairs the manifold and boundary conditions
are such that the objects are quotients (in the sense of definition
\ref{def:quotient-definition}) of a common manifold by two representations
with isomorphic inductions in some supergroup of symmetries.

It turns out, however, that in order to form a quotient by a general
representation we need more singularities than just boundaries (at
least via our construction). A graph is a one dimensional manifold
when all of its vertices are of degree two, and a manifold with boundary
when all vertices are of degree at most two. Unfortunately, even if
a graph has one of these properties, its quotient by a multidimensional
representation (as constructed in section \ref{sec:quotient graph})
need not have either, since the degrees of the vertices are multiplied,
in general, by the dimension of the representation.

Carrying over the construction method of section \ref{sec:quotient graph}
to general Riemannian manifolds (e.g., by replacing graphs with higher
dimensional simplicial structures) yields objects we might call {}``quantum-stratifolds''.
In general, these consist of several Riemannian manifolds of the same
dimension {}``glued'' along their boundaries by homogeneous boundary
conditions (so in dimension one, we obtain the notion of quantum graphs).
When a boundary condition involves the boundaries of more than two
manifolds, the result is no longer a manifold, but rather a stratifold.
Even though this is in general the case, by choosing an appropriate
action, representation and bases, it is possible to obtain manifolds
even when taking a quotient by a multidimensional representation.

\subsection{\label{sub:squareandtriangle}Isospectral drums}

In \cite{JLNP,LPP}, Jakobson et al., and Levitin et al., respectively,
obtain several examples of isospectral domains with mixed Dirichlet-Neumann
boundary conditions, all of which can be interpreted as quotients
with respect to representations sharing a common induction. As a basic
demonstration of the generalization of our theory to higher dimensions,
we reconstruct an isospectral pair consisting of a square and a triangle
with mixed boundary conditions (figure 1 in \cite{LPP}, \ref{fig:square_triangle}
here).

\begin{center}
\begin{figure}[!h]
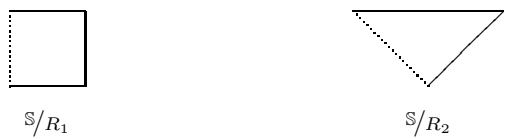

\begin{centering}
\begin{minipage}[c]{0.4\columnwidth}%
\[
\xy(15,15)*{};(15,25)*{}**\dir{.};(15,15)*{};(25,15)*{}**\dir{-};(15,25)*{};(25,25)*{}**\dir{-};(25,15)*{};(25,25)*{}**\dir{-};\endxy\]

\begin{center}
$\nicefrac{\mathbb{S}}{R_{1}}$ 
\par\end{center}%
\end{minipage}%
\begin{minipage}[c]{0.4\columnwidth}%
\[
\xy(5,25)*{};(25,25)*{}**\dir{-};(5,25)*{};(15,15)*{}**\dir{.};(15,15)*{};(25,25)*{}**\dir{-};\endxy\]

\begin{center}
$\nicefrac{\mathbb{S}}{R_{2}}$
\par\end{center}%
\end{minipage}
\par\end{centering}

\caption{The two isospectral domains presented in \cite{LPP}, obtained as
quotients of the square $\mathbb{S}$ (figure \ref{fig:square-and-reflections})
by the representations in \eqref{eq:r1_rep}, \eqref{eq:r2_rep}.
Solid lines indicate Dirichlet boundary conditions and dotted ones
Neumann.}

\label{fig:square_triangle} 
\end{figure}

\par\end{center}

This example rests upon our acquaintance $D_{4}$, so that we can
reuse the definitions and results of section \ref{sec:isospectral_qg}.
In place of the graph in figure \ref{fig:full_square}, we now consider
the full square $\mathbb{S}$, with Dirichlet boundary conditions,
and with $G=D_{4}$ acting as one would expect (figure \ref{fig:square-and-reflections}).

\begin{center}
\begin{figure}[!h]
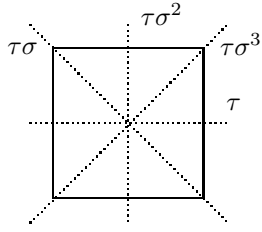

\[
\xy(5,5)*{}="A";(5,25)*{}="B";(25,5)*{}="C";(25,25)*{}="D";"A";"B"**\dir{-};"A";"C"**\dir{-};"C";"D"**\dir{-};"B";"D"**\dir{-};(15,2)*{};(15,28)*{}**\dir{.};(2,15)*{};(28,15)*{}**\dir{.};(19.5,29.5)*{\tau\sigma^{2}};(29,17.5)*{\tau};(2,28)*{};(28,2)*{}**\dir{.};(2,2)*{};(28,28)*{}**\dir{.};(1,25)*{\tau\sigma};(30,25.5)*{\tau\sigma^{3}};\endxy\]

\caption{The square $\mathbb{S}$, and the axes of reflection elements in $D_{4}$.}

\label{fig:square-and-reflections}
\end{figure}

\par\end{center}

The domains in figure \ref{fig:square_triangle} are quotients of
the square $\mathbb{S}$ (figure \ref{fig:square-and-reflections})
by the representations $R_{1}$ and $R_{2}$ of $H_{1},H_{2}\leq G$,
which are defined in \eqref{eq:r1_rep}, \eqref{eq:r2_rep}. Since
$\ind_{H_{1}}^{G}R_{1}\cong\ind_{H_{2}}^{G}R_{2}$, the two domains
are isospectral.

We demonstrate the construction of $\nicefrac{\mathbb{S}}{R_{1}}$.
Recalling that $\mathrm{Hom}_{\mathbb{C}H_{1}}\left(R_{1},\mathcal{H}\left(\mathbb{S}\right)\right)\cong\left(\mathcal{H}\left(\mathbb{S}\right)\right)^{R_{1}}$,
we have again that $\mathcal{H}\left(\nicefrac{\mathbb{S}}{R_{1}}\right)$
should encode the $R_{1}$-isotypic component of $\mathcal{H}\left(\mathbb{S}\right)$.
$\mathbb{T}$, the first quadrant of $\mathbb{S}$ (figure \ref{fig:square_mod_R1}(a)),
is a fundamental domain for the action of $H_{1}$, so that given
$f\in\mathcal{H}\left(\mathbb{T}\right)$ it is possible to construct
at most one function in $\left(\mathcal{H}\left(\mathbb{S}\right)\right)^{R_{1}}$
whose restriction to $\mathbb{T}$ is $f$. Thus, the restriction
map $\Psi:\left(\mathcal{H}\left(\mathbb{S}\right)\right)^{R_{1}}\rightarrow\mathcal{H}\left(\mathbb{T}\right)$
is injective. In order for it to be surjective, we must impose suitable
boundary conditions on $\mathbb{T}$. From \eqref{eq:r1_rep} we obtain
information on $\tilde{f}\in\left(\mathcal{H}\left(\mathbb{S}\right)\right)^{R_{1}}$.
Since such $\tilde{f}$ is anti-symmetric with respect to the action
of $\tau$, it must vanish at the horizontal axis of reflection, and
therefore every $f\in\mathrm{im}\,\Psi$ vanishes at the lower edge
of $\mathbb{T}$. Similarly, every $\tilde{f}\in\left(\mathcal{H}\left(\mathbb{S}\right)\right)^{R_{1}}$
is symmetric with respect to $\tau\sigma^{2}$, so that its normal
derivative at the vertical axis of reflection is zero, and thus all
functions in $\mathrm{im}\,\Psi$ have vanishing normal derivatives
at the left edge of $\mathbb{T}$. This information, summarized in
figure \ref{fig:square_mod_R1}(a), suggests the domain presented
in figure \ref{fig:square_mod_R1}(b) as the quotient $\nicefrac{\mathbb{S}}{R_{1}}$:
a square identical to $\mathbb{T}$, three of whose edges have Dirichlet
boundary condition and one Neumann.

\begin{center}
\begin{figure}[h]
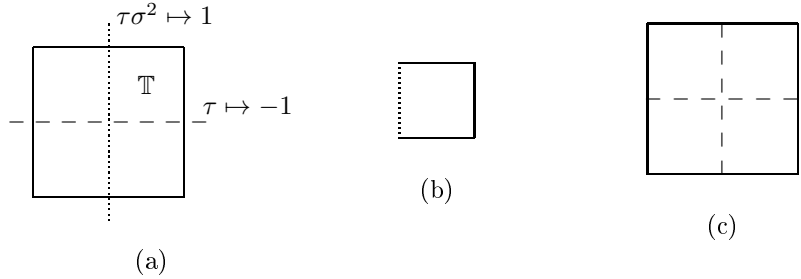

\begin{centering}
\begin{minipage}[c]{0.35\columnwidth}%
\[
\xy(5,5)*{}="A";(5,25)*{}="B";(25,5)*{}="C";(25,25)*{}="D";"A";"B"**\dir{-};"A";"C"**\dir{-};"C";"D"**\dir{-};"B";"D"**\dir{-};(15,2)*{};(15,28)*{}**\dir{.};(2,15)*{};(28,15)*{}**\dir{--};(22.5,29.5)*{\tau\sigma^{2}\mapsto1};(33.5,17)*{\tau\mapsto-1};(20,20)*{\mathbb{T}};\endxy\]

\begin{center}
(a) 
\par\end{center}%
\end{minipage}%
\begin{minipage}[c]{0.25\columnwidth}%
\[
\xy(15,15)*{};(15,25)*{}**\dir{.};(15,15)*{};(25,15)*{}**\dir{-};(15,25)*{};(25,25)*{}**\dir{-};(25,15)*{};(25,25)*{}**\dir{-};\endxy\]

\medskip{}

\begin{center}
(b) 
\par\end{center}%
\end{minipage}%
\begin{minipage}[c]{0.35\columnwidth}%
\begin{center}
$\xy(5,5)*{}="A";(5,25)*{}="B";(25,5)*{}="C";(25,25)*{}="D";"A";"B"**\dir{-};"A";"C"**\dir{-};"C";"D"**\dir{-};"B";"D"**\dir{-};(15,5)*{};(15,25)*{}**\dir{--};(5,15)*{};(25,15)*{}**\dir{--};\endxy$
\par\end{center}

\medskip{}

\begin{center}
(c) 
\par\end{center}%
\end{minipage}
\par\end{centering}

\caption{(a) The fundamental domain $\mathbb{T}$ for $\nicefrac{\mathbb{S}}{H_{1}}$;
every $\tilde{f}\in\left(\mathcal{H}\left(\mathbb{S}\right)\right)^{R_{1}}$
vanishes along the dashed line and has zero normal derivative at the
dotted line. (b) The quotient planar domain $\nicefrac{\mathbb{S}}{R_{1}}$
which encodes this information. The solid lines represent Dirichlet
boundary conditions and the dotted one Neumann. (c) The square $\mathbb{S}'$
of which (b) is a smooth quotient; along the dashed lines functions
need only be once continuously differentiable.}

\label{fig:square_mod_R1} 
\end{figure}

\par\end{center}

Once these boundary conditions are imposed on $\nicefrac{\mathbb{S}}{R_{1}}$,
$\Psi$ is indeed onto: for $f\in\mathcal{H}\left(\nicefrac{\mathbb{S}}{R_{1}}\right)$
which obeys them, we define a function $\tilde{f}$ on $\mathbb{S}$
by $\tilde{f}\Big|_{\mathbb{T}}=f$, $\tau\tilde{f}=-\tilde{f}$,
$\tau\sigma^{2}\tilde{f}=\tilde{f}$, $\sigma^{2}\tilde{f}=-\tilde{f}$.
While $\tilde{f}$ is well defined on the vertical $\tau\sigma^{2}$-axis
even if $f$ does not obey any boundary conditions, it is the requisition
that $f$ vanish on the lower edge of $\mathbb{T}$ that guarantees
that $\tilde{f}$ is well defined on the horizontal $\tau$-axis.
In a similar manner, while at the $\tau$-axis the two one-sided normal
derivatives of $\tilde{f}$ agree a priori, it is the Neumann condition
at the left edge of $\mathbb{T}$ which ensures this at the $\tau\sigma^{2}$-axis.
The boundary conditions thus assure that $\tilde{f}$ is well defined
and continuously differentiable, and being piecewise smooth and a
sum of Laplacian eigenfunctions, it is smooth, and therefore in $\mathcal{H}\left(\mathbb{S}\right)$,
so that $f=\Psi\tilde{f}\in\mathrm{im}\,\Psi$. As $\Psi$ and its
inverse are obviously $\mathbb{C}\left[x\right]$-linear, we have
established $\mathrm{Hom}_{\mathbb{C}H_{1}}\left(R_{1},\mathcal{H}\left(\mathbb{S}\right)\right)\cong\left(\mathcal{H}\left(\mathbb{S}\right)\right)^{R_{1}}\cong\mathcal{H}\left(\nicefrac{\mathbb{S}}{R_{1}}\right)$,
as the definition of a $\nicefrac{\mathbb{S}}{R_{1}}$-domain in \eqref{eq:manifold-quotient-def}
calls for.

Analogously, from the properties of $\tilde{f}\in\left(\mathcal{H}\left(\mathbb{S}\right)\right)^{R_{2}}$
we can deduce the corresponding quotient $\nicefrac{\mathbb{S}}{R_{2}}$.
This process is summarized in the two parts of figure \ref{fig:square_mod_R2}.

\begin{center}
\begin{figure}[!h]
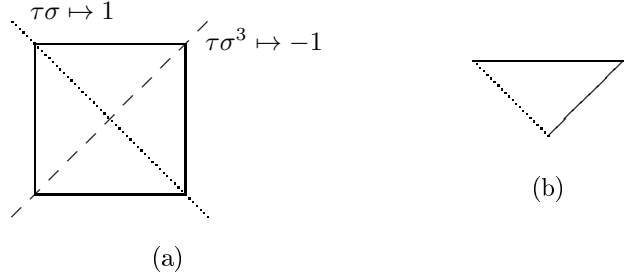

\begin{centering}
\begin{minipage}[c]{0.4\columnwidth}%
\[
\xy(5,5)*{}="A";(5,25)*{}="B";(25,5)*{}="C";(25,25)*{}="D";"A";"B"**\dir{-};"A";"C"**\dir{-};"C";"D"**\dir{-};"B";"D"**\dir{-};(2,28)*{};(28,2)*{}**\dir{.};(2,2)*{};(28,28)*{}**\dir{--};(10,29.5)*{\tau\sigma\mapsto1};(35.5,25.5)*{\tau\sigma^{3}\mapsto-1};\endxy\]

\begin{center}
(a) 
\par\end{center}%
\end{minipage}%
\begin{minipage}[c]{0.4\columnwidth}%
\[
\xy(5,25)*{};(25,25)*{}**\dir{-};(5,25)*{};(15,15)*{}**\dir{.};(15,15)*{};(25,25)*{}**\dir{-};\endxy\]

\medskip{}

\begin{center}
(b) 
\par\end{center}%
\end{minipage}
\par\end{centering}

\caption{(a) The information we have on $\tilde{f}\in\left(\mathcal{H}\left(\mathbb{S}\right)\right)^{R_{2}}$:
it vanishes along the dashed line and has zero normal derivative at
the dotted line. (b) The quotient planar domain $\nicefrac{\mathbb{S}}{R_{2}}$
which encodes this information. The solid lines represent Dirichlet
boundary conditions and the doted one Neumann.}

\label{fig:square_mod_R2}
\end{figure}

\par\end{center}

We return once more to the question of smooth quotients, which was
raised in section \ref{sub:smooth-quotient}. It turns out that even
with the boundary conditions we have imposed on $\nicefrac{\mathbb{S}}{R_{1}}$,
$\Psi$ is not surjective as a function from $\left(C^{^{\infty}}\negthickspace\left(\mathbb{S}\right)\right)^{R_{1}}$
to $C^{^{\infty}}\negthickspace\left(\nicefrac{\mathbb{S}}{R_{1}}\right)$.
Consider for example the smooth function $x^{2}y^{2}$, when regarding
the lower-left corner of $\mathbb{T}$ as the origin%
\footnote{It is not really a function on $\mathbb{T}$, since it does not vanish
at the right and upper edges. This can be rectified by taking $x^{2}y^{2}\left(x^{2}-1\right)\left(y^{2}-1\right)$
instead, but it would clutter the argument.%
}. It is not the restriction of any function in $\left(C^{^{\infty}}\negthickspace\left(\mathbb{S}\right)\right)^{R_{1}}$.
It \textbf{is} the restriction of a function in $\left(C^{1}\!\left(\mathbb{S}\right)\right)^{R_{1}}$,
namely, $x^{2}y\cdot\left|y\right|$. Once again we encounter the
problem of modeling smooth structures by differential boundary conditions
of degree one, which was addressed in section \ref{sec:graphs}. In
fact, $\Psi$ establishes an isomorphism $C^{^{\infty}}\negthickspace\left(\nicefrac{\mathbb{S}}{R_{1}}\right)\cong\left(C^{^{\infty}}\negthickspace\left(\mathbb{S}'\right)\right)^{R_{1}}$,
where $\mathbb{S}'$ is $\mathbb{S}$ after the addition of two {}``Neumann
lines of degree two'' at the axes of reflection corresponding to
$\tau$ and $\tau\sigma^{2}$ (figure \ref{fig:square_mod_R1}(c));
it is a square divided into four, such that a smooth function on $\mathbb{S}'$
is by definition a function which is smooth on each inner (closed)
square, and satisfies the Neumann condition where two squares meet,
or equivalently, is $C^{1}$ at the $\tau$ and $\tau\sigma^{2}$
axes. Had we proceeded by the rigorous method presented in section
\ref{sub:construction-proof}, we would have added these Neumann lines
before the construction, in order to allow a fundamental domain which
is a subcomplex%
\footnote{Each pair of opposite sides in $\mathbb{S}$ is swapped by some element
in $H_{1}$, and as in dimension one we did not allow a vertex to
be moved to a neighbor by a group element, neither should be lines
allowed to, in dimension two. In a more general view, $ $for $\mathrm{Fix}_{\tau}\mathbb{S}$
and $\mathrm{Fix}_{\tau\sigma^{2}}\mathbb{S}$ to be each a subcomplex
of $\mathbb{S}$, $\mathbb{S}$ must be subdivided (by Neumann lines)
into $\mathbb{S}'$.%
}. Again, this could be justified by the preservation of spectral properties:
we have, as for graphs, that $\mathcal{H}\left(\mathbb{S}\right)=\mathcal{H}\left(\mathbb{S}'\right)$,
so that a $\nicefrac{\mathbb{S}}{R}$-quotient (now in the sense of
definition \ref{def:quotient-definition}) is the same thing as a
$\nicefrac{\mathbb{S}'}{R}$-quotient.

We remark that the various constructions demonstrated in section \ref{sec:isospectral_qg}
can be applied analogously to $\mathbb{S}$, enriching the isospectral
pair in figure \ref{fig:square_triangle}. For example, $\nicefrac{\mathbb{S}}{R_{3}}$
would be an orbifold with a line that applies a factor of $i$ to
functions crossing it. The other isospectral families in \cite{JLNP,LPP}
can be obtained from various representations of the general dihedral
groups $D_{n}$, and of the product $D_{4}\times D_{4}$. The interested
reader will find some of these constructions in \cite{BPbS}.

\subsection{\label{sub:drums}The Gordon-Webb-Wolpert drums}

In a similar fashion, we can apply our method to the Gordon-Webb-Wolpert
construction \cite{Gordon1,Gordon2}, obtaining their isospectral
planar domains with new boundary conditions. We follow the exposition
of Buser et al.\ \cite{BuserConway}, who obtain the mentioned drums
as follows: they consider $G_{0}$, a group of motions of the hyperbolic
plane $\mathbb{H}$ ($*444$ in Conway's orbifold notation), and an
epimorphism $\pi:G_{0}\twoheadrightarrow G={\rm PSL}_{3}(2)$. In
$G$ they exhibit two subgroups $A$ and $B$, each isomorphic to
$S_{4}$, that satisfy the Sunada condition \cite{Sunada} with respect
to $G$. The quotients of $\mathbb{H}$ by $\pi^{-1}(A)$ and $\pi^{-1}(B)$
are isometric domains. Both are composed of seven copies of a hyperbolic
triangle (which is a fundamental domain for the action of $G_{0}$),
assembled in different configurations (which are determined by the
coset structure of the pre-images). Finally, by replacing the fundamental
hyperbolic triangle with a suitable Euclidean one, the non-isometric
isospectral drums of Gordon et al.\ are obtained.

An elegant formulation of the Sunada condition for $H_{1}$ and $H_{2}$
in $G$ is that the inductions of the trivial representations $\mathbf{1}_{H_{1}}$
and $\mathbf{1}_{H_{2}}$ to $G$ are isomorphic, i.e. \begin{eqnarray}
\ind_{H_{1}}^{G}\mathbf{1}_{H_{1}} & \cong & \ind_{H_{2}}^{G}\mathbf{1}_{H_{2}}\,.\label{eq:sunada_cond}\end{eqnarray}
 In fact, the connection between $A$ and $B$ is stronger than this
(reflecting a line-point duality in the Fano plane): it turns out
that for every representation $R$ of $S_{4}$, $\ind_{A}^{G}R\cong\ind_{B}^{G}R$.
For each such $R$, we can thus construct an isospectral pair by taking
the quotient of $\mathbb{H}$ by the pullbacks of $R$ to $\pi^{-1}(A)$
and $\pi^{-1}(B)$. Taking $R=\mathbf{1}_{S_{4}}$ will produce once
again the planar drums of Gordon et al. In fact, we shall see in section
\ref{sub:Sunada} that taking quotient (in our sense) by the trivial
representation of a group is equal to taking quotient (in the classical
sense) by the group. Taking $R$ to be the sign representation of
$S_{4}$, and again replacing the fundamental hyperbolic triangles
with Euclidean ones, we obtain the same drums but with different boundary
conditions (figure \ref{fig:new_Gordon}).

\begin{center}
\begin{figure}[!h]
\hfill{}%
\begin{minipage}[t]{0.4\columnwidth}%
\begin{center}
\includegraphics[width=0.7\textwidth]{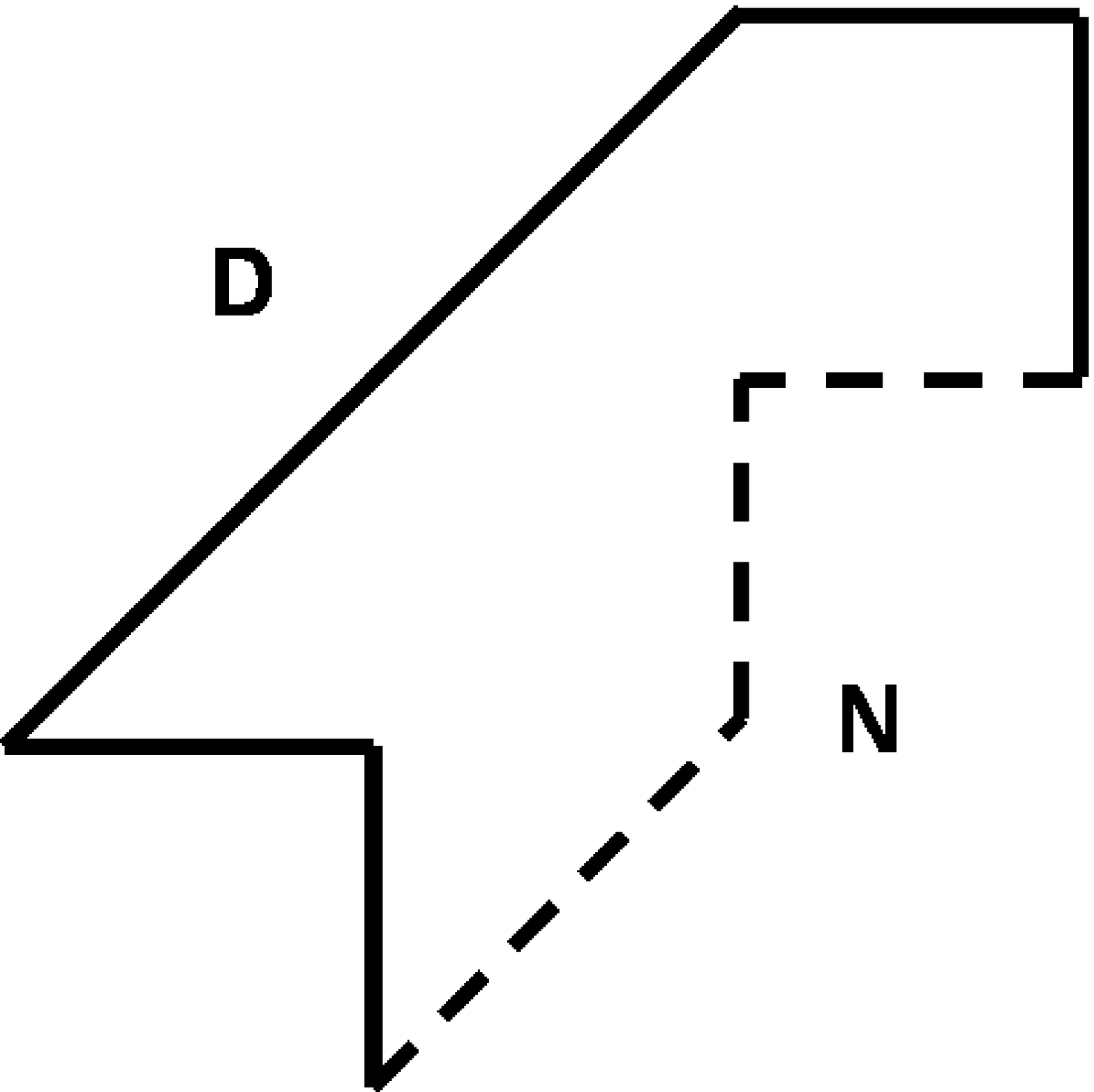} 
\par\end{center}

\begin{center}
(a) 
\par\end{center}%
\end{minipage}\hfill{}%
\begin{minipage}[t]{0.4\columnwidth}%
\begin{center}
\includegraphics[width=0.7\textwidth]{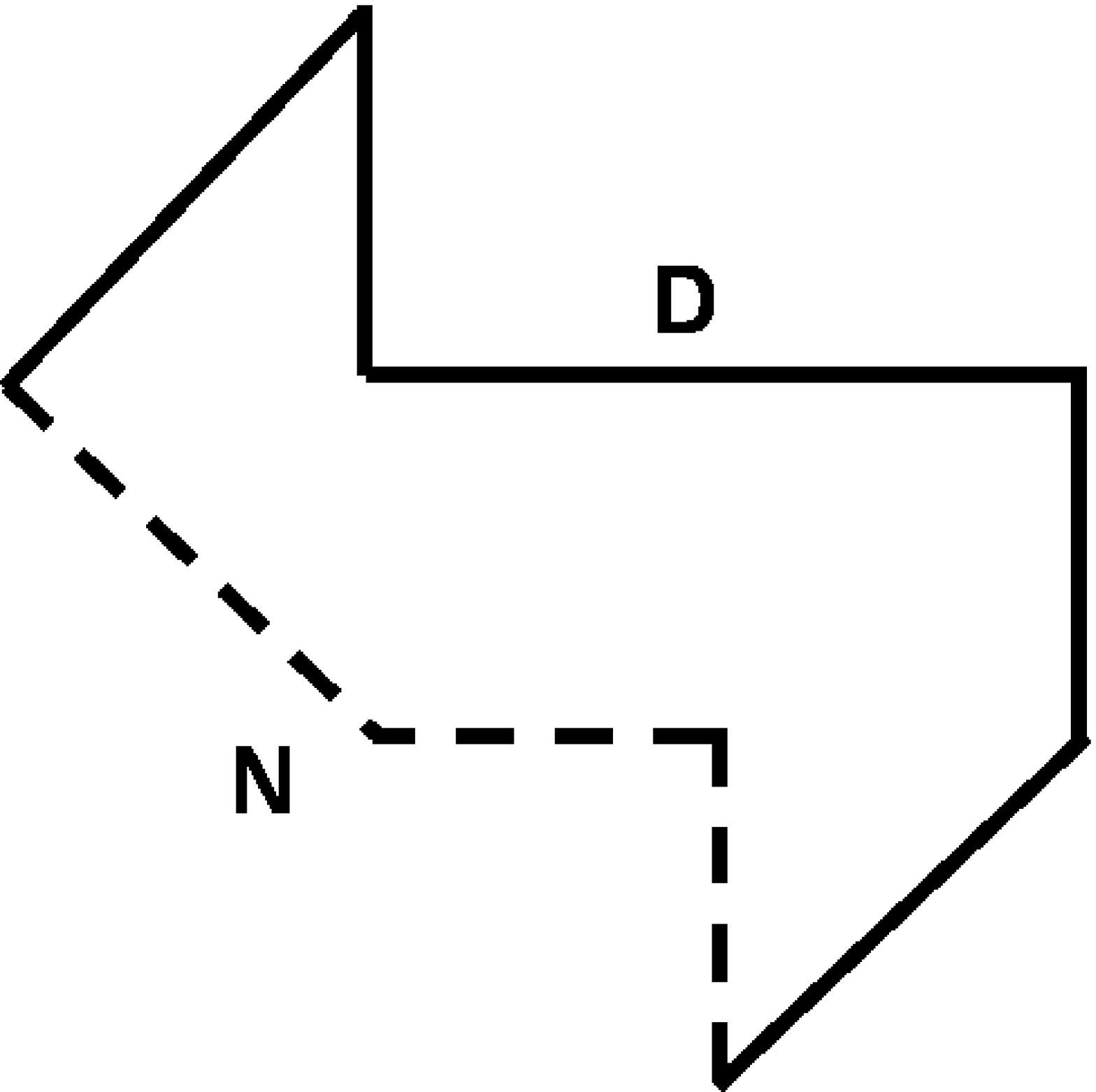} 
\par\end{center}

\begin{center}
(b) 
\par\end{center}%
\end{minipage}\hfill{}

\caption{The isospectral drums of Gordon et al.\ with new boundary conditions.}

\label{fig:new_Gordon} 
\end{figure}

\par\end{center}

We conclude this example by pointing out that in \cite{BuserConway}
a wide variety of isospectral pairs is presented, using various symmetry
groups of $\mathbb{H}$. All these examples can be exploited to construct
other isospectral pairs, as isomorphic inductions may be found either
from Sunada triples or by taking appropriate sums of irreducible representations.

\subsection{\label{sub:Sunada}The Sunada method}

We recall the classical theorem of Sunada \cite{Sunada}:
\begin{quotation}
If $G$ acts freely on a Riemannian manifold $\Gamma$, and $H_{1},H_{2}\leq G$
satisfy \eqref{eq:sunada_cond}, then $\nicefrac{\Gamma}{H_{1}}$
and $\nicefrac{\Gamma}{H_{2}}$ are isospectral manifolds. 
\end{quotation}
Sunada's theorem follows from the definition and corollary at the
beginning of the current section, once we show that for a finite group
$G$ acting freely on a manifold $\Gamma$, the quotient manifold
$\nicefrac{\Gamma}{G}$ is a $\nicefrac{\Gamma}{\mathbf{1}_{G}}$-manifold,
that is,\begin{equation}
\mathcal{H}\left(\nicefrac{\Gamma}{G}\right)\cong\mathrm{Hom}_{\mathbb{C}G}\left(\mathbf{1}_{G},\mathcal{H}\left(\Gamma\right)\right)\,.\label{eq:sunada-our}\end{equation}

This follows from the observation that $\mathrm{Hom}_{\mathbb{C}G}\left(\mathbf{1}_{G},C^{^{\infty}}\negthickspace\left(\Gamma\right)\right)$
corresponds naturally to $C^{^{\infty}}\negthickspace\left(\Gamma\right)^{\mathbf{1}_{G}}=C^{^{\infty}}\negthickspace\left(\Gamma\right)^{G}$,
the trivial component of $C^{^{\infty}}\negthickspace\left(\Gamma\right)$,
and this is the space of functions on $\Gamma$ which are stable under
all elements of $G$. But these are exactly the functions which factor
through $\nicefrac{\Gamma}{G}$, hence $C^{^{\infty}}\negthickspace\left(\nicefrac{\Gamma}{G}\right)\cong\mathrm{Hom}_{\mathbb{C}G}\left(\mathbf{1}_{G},C^{^{\infty}}\negthickspace\left(\Gamma\right)\right)$,
and in particular \eqref{eq:sunada-our} follows.
\begin{rem*}
We can view the preceding argument as yet another proof for Sunada's
theorem, but this would be presumptuous. In fact, Pesce \cite{Pesce}
uses Frobenius Reciprocity in exactly the same manner to reprove Sunada's
theorem. A survey of different proofs for Sunada's theorem, among
them Pesce's, can be found in \cite{Brooks-Sun}.
\end{rem*}

\section{Summary and open questions}

\label{sec:summary}The main construction presented in this paper
is that of objects denoted $\nicefrac{\Gamma}{R}$, where $R$ is
a complex representation of a finite group acting on a geometric object
$\Gamma$. For such $\Gamma$ and $R$ there can be, in general, many
objects so denoted, and they are all isospectral to one another. Furthermore,
these objects are defined so that whenever $\mathrm{Hom}_{\mathbb{C}H_{1}}\left(R_{1},\_\right)\cong\mathrm{Hom}_{\mathbb{C}H_{2}}\left(R_{2},\_\right)$,
where each $R_{i}$ is a representation of a group $H_{i}$ acting
on $\Gamma$, there is also an isospectrality between $\nicefrac{\Gamma}{R_{1}}$
and $\nicefrac{\Gamma}{R_{2}}$. The consequences of this are explored
in section \ref{sec:algebra}, and in particular we find two convenient
means for the construction of isospectral objects:

\textbf{Starting with a group $G$: }take subgroups $H_{1},H_{2}\leq G$
and corresponding representations $R_{1},R_{2}$ sharing a common
induction in $G$~%
\footnote{This resembles the Sunada condition, but is dramatically easier to
achieve, since we are free to take any representations of the subgroups
(instead of only the trivial ones). A systematic approach would be
to take all irreducible characters of subgroups of $G$, induce them
to $G$, find linear dependencies, and sum the corresponding representations
accordingly. Also, any $H_{1}$ and $R_{1}$ are usable with $H_{2}=G$,
by taking $R_{2}=\ind_{H_{1}}^{G}R_{1}$.%
}. For any object $\Gamma$ on which $G$ acts by symmetries, $\nicefrac{\Gamma}{R_{1}}$
and $\nicefrac{\Gamma}{R_{2}}$ are isospectral.

\textbf{Starting with an object $\Gamma$:} find a group $G$ acting
on $\Gamma$ and construct $\nicefrac{\Gamma}{\mathbb{C}G}$ (by some
choice of representatives and bases, as explained in section \ref{sub:construction-proof}).
Any quotient thus obtained is isospectral to $\Gamma$ itself, by
the analogue of proposition \ref{pro:regular-spec} for arbitrary
geometric objects. 

It is natural to ask to what extent the various methods for obtaining
isospectral objects overlap. For example, in section \ref{sec:isospectral_qg},
three isospectral graphs (figure \ref{fig:dihedral_triple}) are obtained
from representations with isomorphic inductions, but at the end of
the same section it is demonstrated that all of them (together with
others) could have also been obtained as $\nicefrac{\Gamma}{R}$ for
a single $R$ (by different choices of bases). Can one expect that
given a basis for $R$, there is always a basis for $\ind_{H}^{G}R$
with respect to which $\nicefrac{\Gamma}{R}$ and $\nicefrac{\Gamma}{{\rm {\rm \mathrm{Ind}}}_{H}^{G}R}$
are isometric?

Even when limiting to the basic quotient construction, questions arise.
For $R$ and $\Gamma$ as above, we have a family of isospectral objects
$\nicefrac{\Gamma}{R}$, varying as one moves between different choices
of bases in the construction, as explained in section \ref{sub:construction-proof}
and demonstrated in the last part of section \ref{sec:isospectral_qg}.
This family has the topology of a manifold, being parametrized by
the action of a general linear group on the space of possible bases.
Surveying this continuum of quotient objects, one might ask where
along it occur changes in the shape of the objects (in contrast with
only boundary condition changes), in the number of connected components,
etc. One can look for certain types of objects in this continuum,
such as manifolds, billiards, objects with real boundary conditions,
or ones with a self-adjoint Laplacian. Such questions seem to lead
to a deeper research in differential and algebraic geometry, investigating
the critical points at which changes occur or the algebraic varieties
at which certain conditions are fulfilled. Except for the basic demonstration
of these phenomena in section \ref{sec:isospectral_qg}, we have not
treated these questions.

We list some more questions that seem interesting, and which we have
not regarded:
\begin{itemize}
\item $\Gamma$ is naturally a $\nicefrac{\Gamma}{\mathbb{C}G}$-graph.
Does it occur by our construction? It seems that the answer is yes,
by taking $G$ as a basis for $\mathbb{C}G$, but we have not shown
this.
\item Given two isospectral objects, can it be decided algorithmically whether
they are representation-quotients of a common object?
\item What are the necessary and sufficient conditions for the quotients
constructed in section \ref{sub:construction-proof} to be proper
quotient graphs (in contrast with generalized ones)? Exact quantum
graphs? Graphs with a self-adjoint Laplacian?
\item Can the isomorphism \eqref{eq:quotient-definition} be made natural,
in a suitable category? This can be interpreted both as (contravariant)
functoriality in $R$, or as functoriality in $\Gamma$, which would
require a definition of quantum graph morphisms.
\item Can the theory presented in this paper be applied to discrete graphs?
To representations of Lie groups acting on Riemannian manifolds?
\item It is clear that $\mathcal{H}\left(\Gamma\coprod\Gamma'\right)=\mathcal{H}\left(\Gamma\right)\oplus\mathcal{H}\left(\Gamma'\right)$,
so that $\sigma_{\Gamma\coprod\Gamma'}\equiv\sigma_{\Gamma}+\sigma_{\Gamma'}$,
and given bases for $R$ and $R'$, their union is a basis for $R\oplus R'$
with respect to which $\nicefrac{\Gamma}{R\oplus R'}$ is isometric
to $\nicefrac{\Gamma}{R}\coprod\nicefrac{\Gamma}{R'}$. Is there an
operation $\otimes$ on graphs, or general geometric object, which
gives $\mathcal{H}\left(\Gamma\otimes\Gamma'\right)=\mathcal{H}\left(\Gamma\right)\otimes\mathcal{H}\left(\Gamma'\right)$,
so that $\sigma_{\Gamma\otimes\Gamma'}\equiv\sigma_{\Gamma}\cdot\sigma_{\Gamma'}$?
What about convolution: $\sigma_{\Gamma\star\Gamma'}\equiv\sigma_{\Gamma}\star\sigma_{\Gamma'}$?
\item A classical conjecture, originally aimed at Riemannian manifolds%
\footnote{This question, in the context of quantum graphs, was suggested to
us by L. Friedlander.%
} \cite{Wig}: for $G={\rm Aut}\Gamma$, and $R=\bigoplus\limits _{i=1}^{r}S_{i}$,
where $\left\{ S_{i}\right\} _{i=1}^{r}$ are the irreducible representations
of $G$, is $\sigma_{\Gamma}^{R}\leq1$? 
\end{itemize}

\section{Acknowledgments}

\label{sec:acknowledgments}It is an honor to acknowledge U. Smilansky,
who is the initiator of this work and an enthusiastic promoter of
it, and a pleasure to thank Z. Sela for his support and encouragement.
We are grateful to M. Sieber for sharing with us his notes, which
led to the construction of the isospectral pair of dihedral graphs.
We are indebted to I. Yaakov whose wise remark has led us to examine
inductions of representations. It is a pleasure to acknowledge G.
Ben-Shach for the fruitful discussions which promoted the research.
We thank D. Schüth for the patient examination of the work and the
generous support. The comments and suggestions offered by J. Brüning,
L. Friedlander, O. Post, Z. Rudnick, and M. Solomyak are highly appreciated.
The work was supported by the Minerva Center for non-linear Physics
and the Einstein (Minerva) Center at the Weizmann Institute, by an
ISF fellowship, and by grants from the GIF (grant I-808-228.14/2003),
and BSF (grant 2006065).

\end{document}